\documentclass[10pt]{amsart}
\usepackage{amsfonts}
\usepackage{color}
\usepackage[square,compress,comma, numbers]{natbib}
\usepackage[colorlinks=true, citecolor=blue, linkcolor=blue]{hyperref}
\allowdisplaybreaks[4]
\usepackage{amssymb}
\usepackage{amsmath}
\definecolor{c20}{rgb}{0.,0.7,0.}
\definecolor{c30}{rgb}{0.,0.,1.}
\definecolor{c40}{rgb}{1,0.1,0.7}
\definecolor{c50}{rgb}{1,0,0}
\definecolor{c60}{rgb}{1,0.9,0.1}
\allowdisplaybreaks[4]

\def\eE#1{\textcolor{c50}{#1}}
\def\eE#1{#1}

\newcommand{\abs}[1]{\left\lvert#1\right\rvert}

\def\cLa#1{\textcolor{c30}{#1}}
\def\cLa#1{#1}

\newcommand{\kb}[1]{\boldsymbol{#1}}
\newcommand{\vk}[1]{\kb{#1}}

\newcommand{\ve}{\varepsilon}

\newcommand{\Abs}[1]{ \biggl \lvert #1 \biggr \rvert}

\newcommand{\E}[1]{\mathbb{E}\left\{#1\right\}}

\newcommand{\pk}[1]{\mathbb{P}\left\{#1 \right\}}

\newcommand{\R}{\mathbb{R}}

\newcommand{\N}{\mathbb{N}}
\newcommand{\inr}{\in \R}

\newcommand{\ldot}{,\ldots,}

\newcommand{\BQN}{\begin{eqnarray}}
\newcommand{\EQN}{\end{eqnarray}}
\newcommand{\BQNY}{\begin{eqnarray*}}
\newcommand{\EQNY}{\end{eqnarray*}}

\newcommand{\BS}{\begin{sat}}
\newcommand{\ES}{\end{sat}}
\newcommand{\BT}{\begin{theo}}
\newcommand{\ET}{\end{theo}}
\newcommand{\BL}{\begin{lem}}
\newcommand{\EL}{\end{lem}}
\newcommand{\BK}{\begin{korr}}
\newcommand{\EK}{\end{korr}}

\newcommand{\BD}{\begin{de}}
\newcommand{\ED}{\end{de}}

\newcommand{\BIT}{\begin{itemize}}
\newcommand{\EIT}{\end{itemize}}
\newcommand{\BDI}{\begin{description}}
\newcommand{\EDI}{\end{description}}

\newcommand{\BRM}{\begin{remarks}}
\newcommand{\ERM}{\end{remarks}}

\newcommand{\BEL}{\begin{lem}}
\newcommand{\EEL}{\end{lem}}

\newtheorem{theo}{Theorem}[section]
\newtheorem{sat}[theo]{Proposition}
\newtheorem{de}[theo]{Definition}
\newtheorem{lem}[theo]{Lemma}

\newtheorem{korr}[theo]{Corollary}
\newtheorem{remark}[theo]{Remark}
\newtheorem{remarks}[theo]{Remarks}
\newtheorem{prop}[theo]{Proposition}

\newcommand{\BPR}{\begin{prop}}
\newcommand{\EPR}{\end{prop}}

\newcommand{\nelem}[1]{{Lemma \ref{#1}}}
\newcommand{\neprop}[1]{{Proposition \ref{#1}}}
\newcommand{\netheo}[1]{{Theorem \ref{#1}}}
\newcommand{\nekorr}[1]{{Corollary \ref{#1}}}

\newcommand{\prooftheo}[1]{ \textsc{\bf Proof of Theorem} \ref{#1}:}
\newcommand{\proofprop}[1]{\textsc{\bf Proof of Proposition} \ref{#1}:}
\newcommand{\prooflem}[1]{\textsc{\bf Proof of Lemma} \ref{#1}:}
\newcommand{\proofkorr}[1]{\textsc{\bf Proof of Corollary} \ref{#1}:}

\newcommand{\COM}[1]{}

\newcommand{\QED}{\hfill $\Box$}

\def\rw{\rightarrow}

\def\IF{\infty}

\def\LT{\left}
\def\RT{\right}

\def\rw{\rightarrow}

\def\X{\vk X}

\def\vn{\varepsilon}
\def\Var{\text{Var}}

\def\Bu+#1{\mathcal{B}^{\varepsilon+}_{u}(#1)}

\def\tla{\tilde{a}}
\def\Q{\mathbb{Q}}
\def\x{\vk{x}}
\def\y{\vk{y}}
\def\z{\vk{z}}

\begin{document}

\title[ Approximation of Kolmogorov-Smirnov Test Statistics]{Approximation of Kolmogorov-Smirnov Test Statistic}

\author{Long Bai}
\address{Long Bai, Department of Actuarial Science, \\
	Faculty of Business and Economics\\
University of Lausanne\\
Chamberonne, 1015 Lausanne, Switzerland
}
\email{Long.Bai@unil.ch}

\author{David Kalaj}
\address{David Kalaj, Faculty of Natural Sciences and Mathematics,
	University of Montenegro, 	Dzordza Vašingtona b.b., Podgorica, Montenegro}
\email{davidk@ac.me}

\bigskip

\date{\today}

 \maketitle
{{\bf Abstract}: Motivated by the  weak limit of the  Kolmogorov-Smirnov test statistics, in this contribution, we concern the asymptotics of
\begin{align*}
\mathbb{P}\left\{\sup_{\boldsymbol{x}\in [0,1]^n}\left(W(\boldsymbol{x})\Big| W(\boldsymbol{1})=w\right)>u\right\}, \ w\in\mathbb{R},
\end{align*}
for large $u$ where $W(\boldsymbol{x})$ is the multivariate Brownian sheet based on a distribution function $F$. The results related to general $F$ are investigated and  some important examples are also showed.


}

\def\indi{\mathbb{I}}
\def\Hn{\mathcal{H}_n}
\section{Introduction}
Let $\X_1 \ldot \X_m$ be independent random vectors with a continuous distribution function (df) $F$. We shall assume without loss of generality that the marginal df's of $F$ are the uniform df on $[0,1]$. Define the empirical df $F_m$ of $F$ by
\BQNY
F_m(\vk{x}):=\frac{1}{m}\sum_{i=1}^m \indi_{[\vk{0},\vk{x}]}(\X_i), \quad \x=(x_1 \ldot x_n) \in [0,1]^n, 
\EQNY
where $\indi_A$ denote the indicator of the set $A$ and $[\vk{0},\vk{x}]=\Pi_{i=1}^n[0,x_i]$. Let below
$W$ denote the (unpinned) Brownian sheet determined by $F$, i.e.,
this is a centered Gaussian random field with covariance function
\BQNY
R(\vk{x},\vk{y})=\E{W(\vk{x})W(\vk{y})}=F(\vk{x}\wedge\vk{y}),\quad
\vk{x},\vk{y}\in [0,1]^n.
\EQNY
Further $W_F$ is the pinned Brownian sheet, which is a centered Gaussian random field  with covariance function
\BQNY
R_F(\vk{x},\vk{y})=\E{W_F(\vk{x})W_F(\vk{y})}=F(\vk{x}\wedge\vk{y})-F(\vk{x})F(\vk{y}),
\quad
\vk{x},\vk{y}\in [0,1]^n.
\EQNY
It is well known, see e.g., \cite{TBSEP1986,JEDKS1952,WCPES1966,MNMS1967} that  $\sqrt{m}(F_m-F)$ converges weakly to $W_F$ as $m\to \IF$  in the space of all bounded functions on $[0,1]^n$ under the topology of uniform convergence. Consequently,  if
\BQNY
T_m^n(F):=\sup_{\vk{x}\in [0,1]^n}\sqrt{m}(F_m(\vk{x})-F(\vk{x}))
\EQNY
is the one-sided Kolmogorov-Smirnov (KS) statistic, then we have the convergence in distribution
\BQNY
T_m^n(F)\rw \sup_{\vk{x}\in [0,1]^n}W_F(\vk{x}),\ m\rw\IF.
\EQNY
For the two-sides KS statistic we have a similar approximation.

For $W_F$, the pinned version of $W$ on $[0,1]^n$, we have the following representation
\BQNY
W_F(\vk{x})=W(\vk{x})-F(\vk{x})W(\vk{1}),\ \vk{x}\in [0,1]^n.
\EQNY
Similarly to Brownian bridge, for $W_F$ we have another conditional representation (see \cite{TBSEP1986}), namely
\BQNY
W_F(\vk{x})=W(\vk{x})\Big| W(\vk{1})=0,\ \vk{x}\in [0,1]^n.
\EQNY
In \cite{ChanLai2006}, as $u\rw\IF$, the asymptotics of
\BQNY
\pk{\sup_{\vk{x}\in [0,1]^n}W_F(\vk{x})>u}=\pk{\sup_{\vk{x}\in [0,1]^n}\LT(W(\vk{x})-F(\vk{x})W(\vk{1})\RT)>u}=\pk{\sup_{\vk{x}\in [0,1]^n}\LT(W(\vk{x})\Big| W(\vk{1})=0\RT)>u}
\EQNY
are studied.
In this paper, we consider a more general case, the asymptotics of
\BQN\label{MMr}
\pk{\sup_{\vk{x}\in [0,1]^n}\LT(W(\vk{x})\Big| W(\vk{1})=w\RT)>u}
\EQN
as $u\rw\IF$ for some constant $w\in\R$.\\ 
Moreover, we give some special cases which can not include in the former senarios.\\
Organization of this paper:  In Section 2 we show our main results and some examples are given in Section 3. Following are the proofs and some useful lemmas in Section 4 and Section 5, respectively.

\COM{Further, if
\BQNY
\widetilde{T}_m^n=\widetilde{T}_m^n(F):=\sup_{\vk{x}\in I^n}\sqrt{m}\abs{F_m(\vk{x})-F(\vk{x})}
\EQNY
is the standard KS statistic, then
\BQNY
\widetilde{T}_m^n(F)\rw \widetilde{M}_F:=\sup_{\vk{x}\in I^n}\abs{W_F(\vk{x})},\ m\rw\IF.
\EQNY
}

\def\seE{ \eE{\mathcal{D}}}
\def\seEd{ \seE_\delta}

\COM{\section{Two dimensional cases}

\BT\label{Thm1}
Suppose that $F(\vk x),x\in [0,1]^2$ is a continuous distribution function. If there \eE{a strictly monotone decreasing} function $h(\eE{x}) \ x\in[S,T]$ for some $0\leq S<T\leq 1$ satisfying
\BQNY
\seE:=\{ \vk x\in [0,1]^2: F(\vk x)=1/2\}=\{\vk{x}=(\eE{x},h(\eE{x})): \forall \eE{x}\in[S,T]\}
\EQNY
and
\BQN \label{Fex}
\lim_{\delta\rw 0}\sup_{\vk{z}\in \seE} \underset{\vk{x}\neq \vk{y}}{\sup_{\abs{\vk{x}-\vk{z}},\abs{\vk{y}-\vk{z}}\leq \delta}}\frac{\abs{ F(\vk{x})-F(\vk{y})-a_1(\vk{z})(x_1-y_1)-a_2(\vk{z})(x_2-y_2)}}
{\abs{x_1-y_1}+\abs{x_2-y_2}}=0,
\EQN
with $a_i(\cdot), i=1,2$  two positive continuous functions, then for any $w\in\R$
\BQN\label{Thm1re1}
\pk{\sup_{\vk{x}\in [0,1]^2}W(\vk{x})>u\Big|W(1,1)=w}\sim
K u^2e^{-2u^2+2uw},  \quad K:= 8\int_{S}^{T} a_1(x,h(x))d x \in (0,\IF).
\EQN
\ET

\BK\label{Corr2}
\eE{If the distribution function $F(\x),\x \in [0,1]^2$ has a bounded positive density function $f$,
\eE{then there exists a strictly monotone decreasing function $h$ on $[S,T] \subset [0,1]$} such that \eqref{Thm1re1} holds with $a_1(x_1,x_2)=\int_{0}^{x_2}f(x_1,t)d t$}.
\EK

\begin{remarks}\label{Rem1}
i) In view of \eqref{Fex} for any $x\in (S,T)$
\BQN
\frac{d h(x)}{d x}=-\frac{a_1(x, h(x))}{a_2(x,h(x))},\ \eE{x\in( S,T)}
\EQN
\eE{implying}
\BQNY
\int_{S}^{T} a_1(x,h(x))d x=
\int_{h(T)}^{h(S)} a_2(\overleftarrow{h}(x),x) d x,
\EQNY
where $\overleftarrow{h}$ is the inverse function of $h$.\\
ii) \cLa{
In \netheo{Thm1}, by the proof, we know that $[0,1]^2$ can be generalised to a manifold $\mathcal{E}\subseteq [0,1]^2$, i.e.
\BQNY
\seE:=\{ \vk x\in \mathcal{E} : F(\vk x)=1/2\}=\{\vk{x}=(\eE{x},h(\eE{x})): \forall \eE{x}\in[S,T]\}.
\EQNY
Then the results are
\BQNY
\pk{\sup_{\vk{x}\in \mathcal{E}}W(\vk{x})>u\Big|W(1,1)=w}\sim
K u^2e^{-2u^2+2uw},  \quad K:= 8\int_{S}^{T} a_1(x,h(x))d x \in (0,\IF).
\EQNY
}
iii)\eE{In \nekorr{Corr2} the assumption of the existence of the density $f$ on $[0,1]^2$ can be weekend to the existence of a positive bounded density on some rectangular set $\mathcal{K}$ that includes  $\seE=\{\x\in [0,1]^2:F(\x)=\frac{1}{2}\}$}.
\end{remarks}}

\def\wx{\widetilde{\vk{x}}}
\def\wy{\widetilde{\vk{y}}}
\def\wz{\widetilde{\vk{z}}}
\section{Main Results}
Before stating our main results, we need to introduce some notation. For $\x,\y\in\R^n$,
\BQNY
&&\vk{x}<\vk{y}\Leftrightarrow x_i<y_i,\quad i=1\ldot   n, \quad \vk{x}\wedge\vk{y}=(x_1\wedge y_1\ldot   x_n\wedge y_n),\\
&&
\vk{x}\pm\vk{y}=(x_1\pm y_1\ldot   x_n\pm y_n),
\quad \abs{\vk{x}}=\sum_{i=1}^n\abs{x_i},\
\x*\y=(x_1\times y_1\ldot x_n\times y_n).
\EQNY
Further, for $\vk{x}<\vk{y}$ we write $[\vk{x},\vk{y}]$ for the set $\Pi_{i=1}^n[x_i,y_i]$ and use $\indi_{A}(\cdot)$ for the indicator function of the set $A\subset [0,1]^n$.\\
Let $\Psi(\cdot)$  denote  the survival function of an $N(0,1)$ random variable.

We write below $\lambda_k(A)$ for the Lebesgue measure on $\R^k$ of some measurable set $A \subset \R^k$.

\BT\label{ThmM1}
Let $F(\vk{x}), \vk x\in [0,1]^n, n\geq 2,$ be a continuous distribution function.  \eE{Suppose that there  exists a  function $h(\wx),\ \wx=(x_1\ldot  x_{n-1}) \in L\subseteq [0,1]^{n-1}$ with $\lambda_{n-1}(L)\neq 0$}
such that
\BQNY
\seE:=\LT\{\vk{x}=(\wx,h(\wx)):\wx\in L\RT\}=\Bigl\{\vk{x}\in [0,1]^n:F(\vk{x})=\frac{1}{2}\Bigr\}
\EQNY
and $\lambda_{n-1}(\partial L)=0$ (i.e., $L$ is a Jordan measurable set).
If $h$ is continuously differentiable in the interior of $L$ and  further
\BQN\label{FEXM}
\lim_{\delta\rw 0}\sup_{\vk{z}\in \seE }\underset{\vk{x}\neq \vk{y}}{\sup_{\abs{\vk{x}-\vk{z}},\abs{\vk{y}-\vk{z}}\leq \delta}}\frac{\abs{F(\vk{x})-F(\vk{y})-\sum_{i=1}^{n}a_i(\vk{z})(x_i-y_i)}}
{\sum_{i=1}^{n}\abs{x_i-y_i}}=0,
\EQN
where $a_i$'s are positive continuous functions, then for $w\in\R$ we have as $u\rw\IF$
\BQN\label{Mulre1}
\pk{\sup_{\vk{x}\in [0,1]^n}W(\vk{x})>u\Big|W(\vk{1})=w}\sim
Ku^{2(n-1)}e^{-2u^2+2uw},
 \EQN
where
\BQN\label{KKt}
\cLa{K:= 2^{3(n-1)}\int_{L} \LT(\Pi_{i=1}^{n-1}a_i(\wx, h(\wx))\RT)d \wx \in \eE{(0,\IF)}}.
\EQN
\ET

\BK\label{Corr3}
If the distribution function $F(\x),\x \in [0,1]^n$ has a bounded positive density function $f$,
then there exists a continuously differentiable function $h$ defined on the interior of a Jordan measurable set $L \subset [0,1]^{n-1}$ with $\lambda_{n-1}(L)>0$ such that \eqref{Mulre1} holds with $a_i(\x)= \partial F(\x) / \partial x_i$.
\EK

Consider the random field $W^*(\vk x) = B_0(F(\vk x)),x\in [0,1]^n$ for a given distribution $F(\vk x), \vk x\in [0,1]^n$ with $B_0$ a standard Brownian bridge. We have that $W^*$ is a centered Gaussian random field, and it has covariance function $F( \vk x) \wedge F(\vk y) -F( \vk x) F(\vk y)$.  In the special case that $F(\vk{x})=F(x_1),\vk{x}\in [0,1]^n$ with \eE{$F$ a univariate distribution}, then $F( \vk x) \wedge F(\vk y)= F( \vk x \wedge \vk y)$, hence $W^*$ has the same law as $W_F$. This observations motivates the result of the next theorem, where essentially we use the fact that the tail asymptotics of the supremum of a Brownian bridge with trend is known, see \cite{GauTrend16}.

\BT\label{Thm2}
Let $F(\x),\x \in [0,1]^n$ be a $n$-dimensional distribution  function. If there exist some $\delta\in (0,\eE{1/2})$ such that
\BQN\label{condition1}
 F(\vk{x}\wedge \vk{y})=F(\vk{x})\wedge F(\vk{y})
\EQN
 holds for $\vk{x}, \vk{y}\in \seEd$ with
\BQNY
\seEd:=\LT\{\vk{x}\eE{\in [0,1]^n}: \frac{1}{2}-\delta \leq F(\vk{x})\leq \frac{1}{2}+\delta \RT\},
\EQNY
then for any $w\inr $
\BQN\label{Thm2re1}
\pk{\sup_{\vk{x}\in [0,1]^n}W(\vk{x})>u\Big|W(\vk{1})=w}\sim e^{-2u^2+2uw}.
\EQN
\ET
\cLa{\begin{remark}
By the proof, we know if \eqref{condition1}  holds for any $\vk{x}, \vk{y}\in [0,1]^n$, then for any $w\in\R$ and $u>0$ we have
\BQNY\label{remeq1}
\pk{\sup_{\vk{x}\in [0,1]^n}W(\vk{x})>u\Big|W(\vk{1})=w}= e^{-2u^2+2uw}.
\EQNY
\end{remark}}
\COM{If \eqref{condition1}  holds for any $\vk{x}, \vk{y}\in [0,1]^n$, then for any $w\in\R$ and $u>0$ we have
\BQN\label{remeq1}
\pk{\sup_{\vk{x}\in [0,1]^n}W(\vk{x})>u\Big|W(\vk{1})=w}= e^{-2u^2+2uw}.
\EQN
and
\BQN\label{remeq2}
\pk{\sup_{\vk{x}\in [0,1]^n}\abs{W(\vk{x})}>u\Big|W(\vk{1})=w}\sim e^{-2u^2+2uw}, \ u\rw\IF.
\EQN
\eE{(14) is exactly  what we have above}
}

Next we give a theorem which further derives the approximation of two-sides KS statistic.
\BT\label{Corr4}
Let $F(\vk{x}), \vk x\in [0,1]^n$ be a continuous distribution function. If further
there exist some $\delta\in (0,1/4)$ such that
\BQN\label{Con2}
 \inf_{\x,\y\in\seEd} F(\x\wedge\y)>\delta
\EQN
holds for $\seEd:=\LT\{\vk{x}\eE{\in \mathcal{E}}: \frac{1}{2}-\delta \leq F(\vk{x})\leq \frac{1}{2}+\delta \RT\}$ where $\mathcal{E}\subseteq[0,1]^n$ satisfies that there exist $\x_0\in\mathcal{E} $ such that $F(\x_0)=\frac{1}{2}$, then we have for $w\in\R$ as $u\rw\IF$
\BQN\label{Mulre2}
&&\pk{\sup_{\vk{x}\in \mathcal{E}}\abs{W(\vk{x})}>u\Big|W(\vk{1})=w}\nonumber\\
&&\quad\quad\sim \pk{\sup_{\vk{x}\in \mathcal{E}}W(\vk{x})>u\Big|W(\vk{1})=w}+\pk{\sup_{\vk{x}\in \mathcal{E}}W(\vk{x})>u\Big|W(\vk{1})=-w}.
\EQN
\ET
\begin{remarks}\label{Rem2}
i) When $\mathcal{E}=[0,1]^n$, for $F(\x)$ satisfying \eqref{condition1} in \netheo{Thm2}, \eqref{Con2} always holds. In fact, for $\x,\y\in\seEd$ with $\delta\in(0,\frac{1}{4})$
\BQNY
F(\x\wedge\y)= F(\x)\wedge F(\y)\geq \frac{1}{2}-\delta>\frac{1}{4}>\delta.
\EQNY
Further by \eqref{Thm2re1} and \eqref{Mulre2}, under the conditions of \netheo{Thm2}, we have as $u\rw\IF$
\BQN\label{Thm2re2}
\pk{\sup_{\vk{x}\in [0,1]^n}\abs{W(\vk{x})}>u\Big|W(\vk{1})=w}\sim c e^{-2u^2+2u\abs{w}},
\EQN
where $c=1$ if $w\neq 0$ and $c=2$ if $w=0$.\\
ii) Under the conditions of \netheo{Thm2}, if further \eqref{Con2} holds for $\mathcal{E}=[0,1]^n$, we have by \eqref{Mulre2} in \netheo{Corr4}
\BQN\label{Mainre2}
\pk{\sup_{\vk{x}\in [0,1]^n}\abs{W(\vk{x})}>u\Big|W(\vk{1})=w}
\sim cKu^{2(n-1)}e^{-2u^2+2uw},
\EQN
where $K$ is the same as in \eqref{KKt} and  $c=1$ if $w\neq 0$ and $c=2$ if $w=0$.
\end{remarks}

\section{Applications}

In this part, we give the asymptotic results of \eqref{MMr} when $F$ are some special cases.
First, we give several two-dimensional cases.
\BPR\label{PROP00}
For $i)-ii)$ below and $w\in \R$, we have that both  \eqref{Mulre1} and \eqref{Mainre2} with $n=2$ hold.\\
i) If $F(\vk{x})=\LT(x_1+x_2-1\RT)^+,\ \vk x \in [0,1]^2$, $K=4$.\\
ii) If $F(\x)=\frac{x_1x_2}{1+(1-x_1)(1-x_2)},\ \vk x \in [0,1]^2$, $K=3\ln 3$.
\EPR
\begin{remarks}
i) In \cite{TBSEP1986}[Theorem 3.1] the following upper bound  for case in \neprop{PROP00} is given
\BQNY
\pk{\sup_{\vk{x}\in [0,1]^2}\LT(W(\vk{x})\Big| W(\vk{1})=0\RT)>u}\leq
\sum_{i=1}^{\IF}(8i^2u^2-2)e^{-2i^2u^2},\quad u>0.
\EQNY
Comparing with our exact result, we see that the prefactor for this upper bound is two times our constant $K=4$.\\
ii) In the light of \cite{TBSEP1986}[Theorem 3.1] for any two-dimensional distribution  $F$ on $[0,1]^2$ and any $u>0$
\BQNY
\pk{\sup_{\vk{x}\in [0,1]^2}W_F(\vk{x})>u}\leq \pk{\sup_{\vk{x}\in [0,1]^2}W_G(\vk{x})>u},
\EQNY
with $G(\vk{x})=(x_1+x_2-1)^{+},\vk x \in [0,1]^2$. Consequently,
our result in \neprop{PROP00} gives an asymptotic upper bound for any $2$-dimensional distribution $F$ on $[0,1]^2$.
\end{remarks}

Following are several multi-dimensional cases.
\BPR\label{PROP0}
For $i)-ii)$ below both  \eqref{Mulre1} and \eqref{Mainre2} hold for any $w\in \R$. Moreover we have: \\
i) If  $F(\vk{x})=\Pi_{i=1}^nx_i, \vk x \in [0,1]^n$, then
 \BQNY
 K=2^{2(n-1)} \int_{L}(\Pi_{i=1}^{n-1}x_i)^{-1}d\wx,
 \EQNY
 and
 \BQNY
 L=\LT\{\wx\in [0,1]^{n-1}:\frac{1}{2}\leq \Pi_{i=1}^{n-1}x_i\leq 1\RT\}.
 \EQNY
 In particular,  $K=4\ln 2$ if $n=2$ and $K=16 (\ln2)^2$ if $n=3$.\\
\COM{ii) If $F(\vk{x})=\frac{1}{n}\sum_{i=1}^nx_i, \vk x \in [0,1]^n$, with
 \BQNY
 \widetilde{K}=2^{3(n-1)}  n^{1-n} \lambda_{n-1}(L),\quad
 L=\LT\{\wx\in [0,1]^{n-1}:\frac{n-2}{2}\leq \sum_{i=1}^{n-1}x_i\leq \frac{n}{2}\RT\}.
 \EQNY
  In particular,  $K=4$ if n=2 and $K=\frac{16}{3} $ if $n=3$. \\}
ii) If
$$F(\x)=d\min_{1\leq i\leq n}x_i+(1-d)\Pi_{i=1}^{n}x_i,\ d\in(0,1),$$
then
\BQNY
K= n(4d)^{(n-1)} \int_{L} \frac{\LT(\Pi_{i=1}^{n-1}x_i\RT)^{n-2}}{\LT(d+(1-d)\Pi_{i=1}^{n-1}x_i\RT)^{n-1}}d \wx,
\EQNY
and
\BQNY
 L=\LT\{\wx\in[0,1]^{n-1}:\frac{1}{2d+2(1-d)\Pi_{i=1}^{n-1}x_i}\leq\min_{1\leq i\leq n-1}x_i\RT\}.
\EQNY
Specially, when $n=2$, we have
\BQNY
K=\frac{8d}{1-d}\ln \LT(\sqrt{1+\frac{1}{(1-d)^2}}-\frac{d}{1-d}\RT).
\EQNY

\EPR

\begin{remark}
i) \eE{The result of \neprop{PROP0}, $i)$ for $n=2$ and $w=0$ agrees with the claim of \cite{LDMSRF1986}[Theorem 1]}.
\COM{\eE{$ii)$  \cite{TBBKI1982} derives the upper and lower bounds for the case
$n=2$ in \neprop{PROP0}, $i)$
Comparing with our results we see that }}
\cLa{ii) In \cite{TBSEP1986}[Theorem 2.1],  a lower bound  for the $n$-dimensional case  is given by
\BQNY
\pk{\sup_{\vk{x}\in [0,1]^n}\LT(W(\vk{x})\Big| W(\vk{1})=w\RT)>u}\geq e^{-2u^2+2uw}
\sum_{i=0}^{n-1}\frac{(2u^2-2uw)^i}{i !},\quad u>w.
\EQNY
Since
$$ \sum_{i=0}^{n-1}\frac{(2u^2-2uw)^i}{i !} \sim \frac{(2u^2)^{n-1}}{(n-1)!}, \quad u\to \IF$$
\eE{comparing with \neprop{PROP0} we  obtain a lower bound for the constant $\Hn $. In the particular case $n=3$ we have $16  (\ln 2)^2 \ge  2$.} \\
}
\end{remark}
\COM{
\BPR\label{KS1}
When $n\geq 2$, we have for $w\in \R$  as $u\rw\IF$
 \eqref{Mulre1} and \eqref{Mulre2} hold with $K=2^{3(n-1)}\Hn \int_{L}n^{1-n}d\wx$ and $L=\LT\{\wx\in I^{n-1}:\frac{n-2}{2}\leq \sum_{i=1}^{n-1}x_i\leq \frac{n}{2}\RT\}$. In particular,  $K=4\ln 2$ if n=2 and $K=16\Hn (\ln2)^2$ if $n=3$.
\EPR

Next we consider the distribution of $M_F$ when $F(\vk{x})=\frac{\sum_{i=1}^nx_i}{n},\ \vk{x}\in I^n$. Let $W_2(\vk{x})$ denote the (unpinned) Brownian sheet on $I^n$, i.e. the zero means Gaussian fields with covariance function
\BQNY
\E{W(\vk{x})W(\vk{y})}= \frac{1}{n}\sum_{i=1}^n(x_i\wedge y_i),\quad \vk{x},\vk{y}\in I^n,
\EQNY
(which means $W(\vk{x})=\frac{1}{\sqrt{n}}\sum_{i=1}^nB^i(x_i)$ where $B^i(x_i)$  are independent Brownian motions.)
and write $W_F$ for the pinned version of $W$ on $I^n$. Then a version of $W_F$ can be obtained from $W$ by the correspondence
\BQNY
W_F(\vk{x})=W(\vk{x})-\LT(\frac{1}{n}\sum_{i=1}^nx_i\RT)W(\vk{1}),\ \ \vk{x}\in I^n.
\EQNY
\BPR\label{KS2}
When $n\geq 2$, we have for $w\in \R$  as $u\rw\IF$
 \eqref{Mulre1} and \eqref{Mulre2} hold with $K=2^{3(n-1)}\Hn \int_{L}n^{1-n}d\wx$ and $L=\LT\{\wx\in I^{n-1}:\frac{n-2}{2}\leq \sum_{i=1}^{n-1}x_i\leq \frac{n}{2}\RT\}$. In particular,  $K=4$ if n=2 and $K=\frac{16}{3}\Hn $ if $n=3$.
\EPR
The lower copula is denoted by $F(x)$, i.e.,
\BQN\label{G}
F(\vk{x})=\LT(\sum_{i=1}^nx_i-1+n\RT)^+, \quad \vk{x}\in I^n.
\EQN
\cite{TBSEP1986} shows for any two-dimensional d.f. $G$ on $I^2$ and for any $u>0$,
\BQNY
\pk{\sup_{\vk{x}\in I^2}W_G(\vk{x})>u}\leq \pk{\sup_{\vk{x}\in I^2}W_F(\vk{x})>u}.
\EQNY
Let $W(\vk{x})$ denote the (unpinned) Brownian sheet on $I^n$, i.e. the zero means Gaussian fields with covariance function
\BQNY
\E{W(\vk{x})W(\vk{y})}= \LT(\sum_{i=1}^nx_i\wedge y_i-1+n\RT)^+,\quad \vk{x},\vk{y}\in I^n,
\EQNY
and write $W_F$ for the pinned version of $W$ on $I^n$. Then a version of $W_F$ can be obtained from $W$ by the correspondence
\BQNY
W_F(\vk{x})=W(\vk{x})-\LT(\sum_{i=1}^nx_i-1+n\RT)^+W(\vk{1}),\ \ \vk{x}\in I^n.
\EQNY
\BPR\label{KS3} When $n\geq 2$, we have for $w\in \R$  as $u\rw\IF$
 \eqref{Mulre1} and \eqref{Mulre2} hold with $K=2^{3(n-1)}\Hn \int_{L}d\wx$ and $L=\LT\{\wx\in I^{n-1}:n-\frac{3}{2}\leq \sum_{i=1}^{n-1}x_i\leq n-\frac{1}{2}\RT\}$. Especially, $K=4$ if n=2 and $K=8\Hn $ if $n=3$.
\EPR
Finially we consider the distribution of $M_F$ when $F(\vk{x})=\min_{1\leq i\leq n}x_i,\ \vk{x}\in I^n$. Let $W(\vk{x})$ denote the (unpinned) Brownian sheet on $I^n$, i.e. the zero means Gaussian fields with covariance function
\BQNY
\E{W(\vk{x})W(\vk{y})}= \min_{1\leq i\leq n}(x_i\wedge y_i),\quad \vk{x},\vk{y}\in I^n.
\EQNY
and write $W_F$ for the pinned version of $W$ on $I^n$. Then a version of $W_F$ can be obtained from $W$ by the correspondence
\BQNY
W_F(\vk{x})=W(\vk{x})-\min_{1\leq i\leq n}(x_i\wedge y_i)W(\vk{1}),\ \ \vk{x}\in I^n.
\EQNY
\BPR\label{KS4}
We have that for $w\in\R$, as $u\rw\IF$ \eqref{Thm2re1} and \eqref{Thm2re2} hold.
\EPR}
Next proposition is a case which satisfies the \eqref{condition1} in \netheo{Thm2}.
\BPR\label{PROP1}
For  $d\in (0,1)$ and
\BQNY
F(\x)=
\left\{
\begin{array}{ll}
\frac{1}{2d}\min_{1\leq i\leq n}x_i,&\ \text{if}\  \min_{1\leq i\leq n}x_i\leq d,\\
\frac{1}{2(1-d)}\min_{1\leq i\leq n}x_i+\frac{1-2d}{2(1-d)},&\ \text{if}\ \min_{1\leq i\leq n}x_i\geq d,
\end{array}
\right.
\EQNY
we have that both  \eqref{Thm2re1} and \eqref{Thm2re2} hold
for any $w\in \R$, i.e.
\BQNY
\pk{\sup_{\vk{x}\in [0,1]^n}\LT(W(\vk{x})\Big|W(\vk{1})=w\RT)>u}= e^{-2u^2+2uw}, \ u>0,
\EQNY
and
\BQNY
\pk{\sup_{\vk{x}\in [0,1]^n}\LT(\abs{W(\vk{x})}\Big|W(\vk{1})=w\RT)>u}\sim \cLa{c}e^{-2u^2+2u\abs{w}}, \ u\rw\IF,
\EQNY
with $c=1$ if $w\neq0$ and $c=2$ if $w=0$.
\COM{iii) For $d\in(0,1)$ and
$$F(\x)=d\min_{1\leq i\leq n}x_i+\frac{(1-d)}{n}\sum_{i=1}^{n}x_i,$$
we have that  \eqref{Mulre1} holds for any $w\in \R$ with
\BQNY
\widetilde{K}= \frac{(8(1-d))^{n-1}}{n^{n-2}}\Hn \lambda_{n-1}(L),
\text{and}\ L=\LT\{\wx\in[0,1]^{n-1}:\frac{n-2(1-d)\sum_{i=1}^{n-1}x_i}{2+2(n-1)d}\leq\min_{1\leq i\leq n-1}x_i\RT\}.
\EQNY
Specially, when $n=2$, we have that \eqref{Thm1re1} holds for any $w\in \R$ with
\BQNY
K=4(1-d).
\EQNY}
\EPR
\begin{remark}
In \neprop{PROP1}, if $d=\frac{1}{2}$, $F(\x)=\min_{1\leq i\leq n}x_i$, i.e., the upper copula.
\end{remark}

\section{Proofs}

\prooftheo{ThmM1} Hereafter, we denote by $\mathbb{Q}_i,\ i\in\mathbb{N}$ some positive constants that may differ from line to line.\\
By the monotonicity and continuity of the distribution function $F(\vk{x}), \vk{x}\in [0,1]^n$, $h(\wx)$ is a continuous function over  $L\subset[0,1]^{n-1}$ which is strictly decrease along every line parallel to the  axes (and so on all increasing paths).\\
Then for $\wx\in L$ and $x_n=h(\wx)$ we set
\BQNY
\tla_i(\wx):=a_i(\wx,h(\wx)),\  i=1\ldot n,
\EQNY
and
 \BQNY
\underline{a}_i:=\inf_{\x\in \seE}a_i(\x)>0,\ \overline{a}_i:=\sup_{\vk{x}\in \seE}a_i(\vk{x})<\IF,  i=1\ldot n,
 \EQNY
where we use the fact that $a_i(\x)$'s are continuous and positive function.

We have for $u>w$
\begin{align*}
\pk{\sup_{\vk{x}\in [0,1]^n}W(\vk{x})>u\Big|W(\vk{1})=w}
&=\pk{\sup_{\vk{x}\in [0,1]^n}(W_F(\vk{x})+F(\vk{x})w)>u},
\end{align*}
where $W_F(\vk{x}):=W(\vk{x})-F(\vk{x})W(\vk{1})$.
The variance function of $W_F(\vk{x})$ is
\BQNY
\sigma^2_F(\vk{x}):=\Var(W_F(\vk{x}))=F(\vk{x})(1-F(\vk{x})),\quad \vk{x}\in [0,1]^n,
\EQNY
which attains its maximum equal to $\frac{1}{4}$ over $[0,1]^n$ at $\vk{z}$ with $F(\vk{z})=\frac{1}{2}$, i.e., at $\seE$ and as $F(\x)\rw \frac{1}{2}$
\BQN\label{var11}
\frac{1}{2}-\sigma_F(\x)\sim\LT(F(\x)-\frac{1}{2}\RT)^2.
\EQN
\COM{Further, by \eqref{FEXM1}, we have that
\BQN\label{hhB}
\lim_{\delta\rw 0}\sup_{\z\in\seE}\sup_{\abs{\x-\z}<\delta}\abs{\frac{\frac{1}{2}-\sigma_F(\vk{x})}
{\LT(\sum_{i=1}^na_i(\vk{z})(x_i-z_i)\RT)^2}-1}=0.
\EQN}
Since by \eqref{FEXM}, there exist $\vn_1\in\LT(0,\min_{1\leq i\leq n}\underline{a}_i\RT),$ for any $\vk{z}\in \seE$, if $\abs{\x-\z},\abs{\y-\z}<\delta$
\BQN
F(\x)-F(\y)\leq \sum_{i=1}^{n}a_i(\z)(x_i-y_i)+\vn_1\sum_{i=1}^n\abs{x_i-y_i}\leq \Q_1\abs{\x-\y},\label{ThmMbo1}\\
F(\x)-F(\y)\geq \sum_{i=1}^{n}a_i(\z)(x_i-y_i)-\vn_1\sum_{i=1}^n\abs{x_i-y_i}\geq \Q_2\abs{\x-\y}\label{ThmMbo2}.
\EQN
Thus $F(\x)$ is strictly increasing along every line parallel to the  axes in
\BQN\label{ED}
E(\delta):=\{\x\in [0,1]^n:\abs{\x-\z}\leq \delta, \z\in \seE\}\supseteq \seE,
\EQN
and for any $\delta_0\in (0,1/2)$ we can take $\delta\in(0,\frac{1}{2})$ small enough such that
\BQN\label{Fb2}
\sup_{\x\in E(\delta)}\abs{F(\vk{x})-1/2}\leq \delta_0,
\EQN
with $\delta_0\rw 0$ as $\delta\rw 0$ and
\BQNY
\LT(\frac{1}{2}-\delta_0\RT)^2\leq \sigma^2_F(\vk{x})\leq 1,\ \x\in E(\delta).
\EQNY
For the correlation function $r_F(\vk{x},\vk{y}):=Cov\LT(\frac{W_F(\vk{x})}{\sigma_F(\vk{x})},
\frac{W_F(\vk{y})}{\sigma_F(\vk{y})}\RT)$, we have for any $\z\in\seE$, if $\abs{\x-\z},\abs{\y-\z}\leq\delta$
\BQN
1-r_F(\vk{x},\vk{y})
&=&1-\frac{\E{W_F(\vk{x})W_F(\vk{y})}}{\sigma_F(\vk{x})\sigma_F(\vk{y})}\label{rrr1}\\
&=&\frac{\E{\LT(W_F(\vk{x})-W_F(\vk{y})\RT)^2}
-\LT(\sigma_F(\vk{x})-\sigma_F(\vk{y})\RT)^2}
{2\sigma_F(\vk{x})\sigma_F(\vk{y})}.\nonumber
\EQN
By \eqref{Fb2}, we have
\BQNY
\LT(\sigma_F(\vk{x})-\sigma_F(\vk{y})\RT)^2
&=&\frac{\LT(\sigma^2_F(\vk{x})-\sigma^2_F(\vk{y})\RT)^2}
{\LT(\sigma_F(\vk{x})+\sigma_F(\vk{y})\RT)^2}\\
&\leq&\frac{1}{(1-2\delta_0)^2}\LT((F(\vk{x})-F(\vk{y}))-(F^2(\vk{x})-F^2(\vk{y}))\RT)^2\\
&=&\frac{1}{(1-2\delta_0)^2}(F(\vk{x})-F(\vk{y}))^2(1-(F(\vk{x})+F(\vk{y})))^2\\
&\leq&\frac{4\delta_0^2\Q_2^2}{(1-2\delta_0)^2}\abs{\x-\y}^2,
\EQNY
and
\BQNY
\E{\LT(W_F(\vk{x})-W_F(\vk{y})\RT)^2}
&=&\E{\LT(W(\vk{x})-W(\vk{y})-(F(\vk{x})-F(\vk{y}))W(1,1)\RT)^2}\\
&=&F(\vk{x})+F(\vk{y})-2F(\vk{x}\wedge \vk{y})-(F(\vk{x})-F(\vk{y}))^2.
\EQNY
Since by \eqref{FEXM}, there exist $\vn_1\in(0,\min_{1\leq i\leq n}\underline{a}_i),$ for any $\vk{z}\in \seE$, if $\abs{\x-\z},\abs{\y-\z}<\delta$
\BQNY
F(\vk{x})+F(\vk{y})-2F(\vk{x}\wedge \vk{y})\leq
\sum_{i=1}^{n}(\tla_i(\wz)+\vn_1)\abs{x_i-y_i},\\
F(\vk{x})+F(\vk{y})-2F(\vk{x}\wedge \vk{y})\geq
\sum_{i=1}^{n}(\tla_i(\wz)-\vn_1)\abs{x_i-y_i}.
\EQNY
\eE{Consequently, for any $\x,\y\in E(\delta) $}
\BQNY
1-r_F(\x,\y)&\leq&
 \frac{2}{(1-2\delta_0)^2}\LT(F(\vk{x})+F(\vk{y})-2F(\vk{x}\wedge \vk{y})\RT)\\
&\leq& \frac{2}{(1-2\delta_0)^2}\LT(\sum_{i=1}^n(\tla_i(\wz)+\vn_1)\abs{x_i-y_i}
\RT)\\
&\leq&2(1+\vn_2)\LT(\sum_{i=1}^n(\tla_i(\wz)+\vn_1)\abs{x_i-y_i}
\RT)
\EQNY
and
\BQNY
1-r_F(\x,\y)&\geq&
 \frac{2}{(1+2\delta_0)^2}\E{\LT(W_F(\vk{x})-W_F(\vk{y})\RT)^2}
 -\frac{4\delta_0^2\Q_2^2}{(1-2\delta_0)^2}\abs{\x-\y}^2\\
 &\geq&
 \frac{2}{(1+2\delta_0)^2}\LT(F(\vk{x})+F(\vk{y})-2F(\vk{x}\wedge \vk{y})-\Q_2^2\abs{\x-\y}^2\RT)
 -\frac{4\delta_0^2\Q_2^2}{(1-2\delta_0)^2}\abs{\x-\y}^2\\
 &\geq&
 \frac{2}{(1+2\delta_0)^2}\LT(\sum_{i=1}^n\tla_i(\wz)\abs{x_i-y_i}\RT)
 -\LT(\Q_2^2+\frac{4\delta_0^2\Q_2^2}{(1-2\delta_0)^2}\RT)\abs{\x-\y}^2\\
 &\geq&
 2(1+\vn_4)\LT(\sum_{i=1}^n\tla_i(\wz)\abs{x_i-y_i}\RT),
\EQNY
where we use the fact that for $\x,\y\in E(\delta)$
\BQNY
\frac{2}{(1+2\delta)^2}\leq \frac{1}{2\sigma_F(\vk{x})\sigma_F(\vk{y})}\leq\frac{2}{(1-2\delta)^2}.
\EQNY
Hence
\BQN\label{r11}
\lim_{\delta\rw 0}\sup_{\vk{z}\in \seE}\underset{\x\neq\y}{\sup_{\abs{\x-\z},\abs{\y-\z}<\delta}} \abs{\frac{1-r_F(\vk{x},\vk{y})}{2\sum_{i=1}^na_i(\vk{z})\abs{x_i-y_i}}-1}=0.
\EQN

Since for $x,y\in(0,1)$
\BQNY
\sqrt{(x-x^2)(y-y^2)}\geq (x\wedge y-xy),
\EQNY
where the \eE{equality} holds only when $x=y$, then for $\vk{x}, \vk{y}\in E(\delta)$
\BQN\label{eqsig}
\sigma_F(\vk{x})\sigma_F(\vk{y})&\geq& F(\vk{x})\wedge F(\vk{y})-F(\vk{x})F(\vk{y})\nonumber\\
&\geq&F(\vk{x}\wedge \vk{y})-F(\vk{x})F(\vk{y})\nonumber\\
&=&\E{W_F(\vk{x})W_F(\vk{y})},
\EQN
and for $\vk{x}\neq \vk{y}$, $\vk{x},\vk{y}\in E(\delta)$, if $F(\vk{x})=F(\vk{y})$, $F(\vk{x}\wedge \vk{y})<F(\vk{x})\wedge F(\vk{y})$ \eE{since $F(\vk{x})$ is strictly increasing along every line parallel to the axes in $E(\delta)$}. Then in \eqref{eqsig},
at least one of the two inequality strictly holds implying
\BQN\label{Thm31ru}
r_F(\vk{x},\vk{y})< 1
\EQN
holds for $\vk{x}, \vk{y}\in E(\delta)$ and $\vk{x}\neq \vk{y}$.\\
\COM{By \eqref{ThmMbo1} and \eqref{ThmMbo2}, for $\wx,\ \wy\in L$ with $x_i=y_i, i=2\ldot  n-1$,
\BQNY
&&0=F(\wx,h(\wx))-F(\wy,h(\wy))\leq \tla_1(\wx)(x_1-y_1)+\tla_n(\wx)(h(\wx)-h(\wy))+\vn_1\abs{x_1-y_1}+\vn\abs{h(\wx)-h(\wy)},\\
&&0=F(\wx,h(\wx))-F(\wy,h(\wy))\geq \tla_1(\wx)(x_1-y_1)+\tla_n(\wx)(h(\wx)-h(\wy))-\vn_1\abs{x_1-y_1}-\vn\abs{h(\wx)-h(\wy)},
\EQNY
i.e.
\BQNY
&&-\tla_n(\wx)(h(\wx)-h(\wy))-\vn\abs{h(\wx)-h(\wy)}\leq\tla_1(\wx)(x_1-y_1)+\vn_1\abs{x_1-y_1},\\
&&-\tla_n(\wx)(h(\wx)-h(\wy))+\vn\abs{h(\wx)-h(\wy)}\geq\tla_1(\wx)(x_1-y_1)-\vn_1\abs{x_1-y_1}.
\EQNY
Then if we set $x_1>y_1$, the right hand side of the last inequality is great than zero and $h(\wx)<h(\wy)$ and
\BQN\label{hhB}
\frac{\tla_1(\wx)-\vn_1}{\tla_n(\wx)+\vn_1}\abs{x_1-y_1} \leq \abs{h(\wx)-h(\wy)}\leq \frac{\tla_1(\wx)+\vn_1}{\tla_n(\wx)-\vn_1}\abs{x_1-y_1}.
\EQN
Thus $h(\wx), \wx\in L^o$ with $L^o$ the interior of $L$ is strictly decreasing continuous differentiable about $x_1$ on $L^o$ with
\BQNY
\frac{\partial h(\wx)}{\partial x_1}=-\frac{\tla_1(\wx)}{\tla_n(\wx)}.
\EQNY
Further, $h(\wx), \wx\in L^o$ is strictly decreasing continuous differentiable on $L^o$ with
\BQN
\frac{\partial h(\wx)}{\partial x_i}=-\frac{\tla_i(\wx)}{\tla_n(\wx)},\ i=1\ldot   n-1.
\EQN}

 Set for $\x\in [0,1]^n$ and $A,\ B\subseteq [0,1]^n$
 $$\rho(\x,A)=\inf_{\y\in A}\abs{\x-\y},\ \rho(A, B)=\inf_{\x\in A, \y\in B}\abs{\x-\y} $$
and
\BQNY
E_0(\delta)=\{\vk{x}:\rho(\wx, L)\leq \delta, \abs{x_n-q(\wx)}\leq \delta\}
\EQNY
where
\BQNY
q(\wx)=\LT\{
\begin{array}{ll}
h(\wx),\ \text{if} \  \wx\in L,\\
h(\wy)+\sum_{i=1}^{n-1}\tla(\wy)(x_i-y_i) , \ \text{if} \  \wx\notin L\ \text{and}\ \wy\in\{\wz:\abs{\wz-\wx}=\rho(\wx,L)\}.
\end{array}
\RT.
\EQNY
Then $E_0(\delta)\supset \seE$.
Since $\sigma^2_F(\vk{x})$ is a continuous function, we have
\BQNY
\sigma^2_m=\sup_{\vk{x}\in [0,1]^n\setminus E_0(\delta)}\sigma^2_F(\vk{x})< \frac{1}{4}.
\EQNY
By again Borell-TIS inequality (ref.\cite{AdlerTaylor}), as $u\rw\IF$
\BQN\label{Thm31p2}
\pk{\sup_{\vk{x}\in [0,1]^n\setminus E_0(\delta) }\LT(W_F(\vk{x})+F(\vk{x})w\RT)>u}
&\leq&\pk{\sup_{\vk{x}\in [0,1]^n\setminus E_0(\delta)}W_F(\vk{x})>u-\abs{w}}\nonumber\\
&\leq&\exp\LT(-\frac{(u-\abs{w}-\mathbb{Q}_2)^2}{2\sigma^2_m}\RT)\nonumber\\
&=&o\LT(u^{2(n-1)}e^{-2u^2+2uw}\RT),
\EQN
where $\mathbb{Q}_2=\E{\sup_{\x\in[0,1]^n}W_F(\x)}<\IF$.\\
We have
\begin{align}
\pk{\sup_{\vk{x}\in [0,1]^n}(W_F(\vk{x})+F(\vk{x})w)>u}
&\leq
\pk{\sup_{\vk{x}\in [0,1]^n\setminus E_0(\delta)}(W_F(\vk{x})+F(\vk{x})w)>u} +\Pi_1(u),\label{Thm31ub1}
\end{align}
and
\begin{align}
\pk{\sup_{\vk{x}\in [0,1]^n}(W_F(\vk{x})+F(\vk{x})w)>u}
&\geq\Pi_1(u),\label{Thm31lb2}
\end{align}
where
\BQNY
\Pi_1(u)=\pk{\sup_{\vk{x}\in E_0(\delta)}(W_F(\vk{x})+F(\vk{x})w)>u}.
\EQNY

Next we consider $\Pi_1(u)$.
We have
\BQNY
\Pi_1(u)=\pk{\sup_{\vk{x}\in E_0(\delta)}(X(\vk{x})+g(\vk{x}))>\mu}
\EQNY
where $\mu=2u-w,\ \  g(\vk{x})=w(2F(\x)-1)$
 and $X(\x)=2W_F(\x).$

\def\wz{\widetilde{\vk{z}}}

Notice that the variance function $\sigma^2_X(\vk{x})$ of $X(\vk{x})$ attains its maximum equal to $1$ at $\seE$. In \eqref{FEXM} for $\x\in E_0(\delta)$ taking $\y=\z=(\wx,h(\wx))$ leading
\BQN\label{con1}
\lim_{\delta\rw 0}\sup_{\wx\in L }\underset{x_n\neq h(\wx)}{\sup_{\abs{x_n-h(\wx)}\leq \delta}}\frac{\abs{F(\vk{x})-1/2-\tla_n(\wx)\LT(x_n-h(\wx)\RT)}}
{\abs{x_n-h(\wx)}}=0,
\EQN
which combined with \eqref{var11} implies
\BQN\label{Thm31sig}
\lim_{\delta\rw 0}\sup_{\wx\in L}\underset{x_n\neq h(\wx)}{\sup_{\abs{x_n-h(\wx)}<\delta}}\abs{\frac{1-\sigma_X(\vk{x})}
{2\tla_n^2(\wx)(x_n-h(\wx))^2}-1}=0
\EQN
and
\BQN\label{Thm31gg}
\lim_{\delta\rw0}\sup_{\wx\in L}\underset{x_n\neq h(\wx)}{\sup_{\abs{x_n-h(\wx)}<\delta}}
\abs{\frac{g(\vk{x})}{2w\tla_n(\wx)(x_n-h(\wx))}-1}=0.
\EQN
Further in \eqref{FEXM} for $\x,\y\in \seE$ taking $\z=\x$  leads
\BQN\label{con1}
\lim_{\delta\rw 0}\sup_{\x,\y\in \seE}\underset{\x\neq\y}{\sup_{\abs{\x-\y}\leq \delta}}\frac{\abs{\sum_{i=1}^{n}\tla_i(\wx)\LT(x_i-y_i\RT)}}
{\sum_{i=1}^n\abs{x_i-y_i}}=0,
\EQN
which derive that for any small $\delta\in(0,1)$ there exist a constant $\mathbb{Q}_3$ such that
\BQNY
\sup_{\wx,\wy\in L}\underset{\wx\neq\wy}{\sup_{\abs{\wx-\wy}\leq \delta}}\frac{\abs{h(\wx)-h(\wy)}}{\sum_{i=1}^{n-1} \abs{x_i-y_i}}
\leq  \mathbb{Q}_3.
\EQNY
Thus for small $\delta\in(0,1)$ there exist a constant $ \mathbb{Q}_4$ such that
\BQN
\sup_{\x,\y\in E_0(\delta)}\underset{\x\neq\y}{\sup_{\abs{\x-\y}\leq \delta}}\frac{\abs{(x_n-h(\wx))-(y_n-h(\wy))}}{\sum_{i=1}^{n} \abs{x_i-y_i}}
\leq  \mathbb{Q}_4.
\EQN
\COM{Recall that
\BQNY
\sup_{\x\in E_0(\delta)}\abs{F(\x)-1/2}\leq \delta_0,
\EQNY
which together with \eqref{ThmMbo1} and \eqref{con1} derive that for $\x,\y\in E_0(\delta)$
\BQN
\abs{\frac{1}{\sigma_X(\x)}-\frac{1}{\sigma_X(\y)}}
&=&\frac{1}{2}\LT(\frac{\sqrt{F(\y)(1-F(\y))}-\sqrt{F(\x)(1-F(\x))}}
{\sqrt{F(\x)F(\y)(1-F(\x))(1-F(\y))}}\RT)\nonumber\\
&\leq&\frac{1}{2}\LT(\frac{F(\y)(1-F(\y))-F(\x)(1-F(\x))}
{2(1-2\delta_0)^3}\RT)\nonumber\\
&\leq&\frac{1}{2}\LT(\frac{(F(\y)-F(\x))(\abs{1/2-F(\x)}+\abs{1/2-F(\y)})}
{2(1-2\delta_0)^3}\RT)\nonumber\\
&\leq&\mathbb{Q}_3 \LT(\abs{x_n-h(\wx)}+\abs{y_n-h(\wy)}\RT)\LT(\sum_{i=1}^n\abs{x_i-y_i}\RT),
\EQN
and
\BQN
\abs{g(\x)-g(\y)}=2\abs{w}\abs{F(\x)-F(\y)}\leq \mathbb{Q}_4\LT(\sum_{i=1}^n\abs{x_i-y_i}\RT).
\EQN}
For the correlation function $r_X(\x,\y)$ of $X(\x)$, by \eqref{r11}
\BQN\label{Thm31r2}
\lim_{\delta\rw 0}\sup_{\z\in \seE}\underset{\x\neq\y}{\sup_{\abs{\x-\z},\abs{\y-\z}<\delta}}
\abs{\frac{1-r_X(\vk{x},\vk{y})}{2\sum_{i=1}^{n}a_i(\z)\abs{x_i-y_i}
}-1}=0
\EQN
and further for any $\abs{\x-\z},\abs{\y-\z}<\delta$ with $\z\in \seE$
\BQN\label{Thm31r3}
2(1-\vn)\sum_{i=1}^{n}a_i(\z)\abs{x_i-y_i}\leq 1-r_X(\vk{x},\vk{y})\leq 2(1+\vn)\sum_{i=1}^{n}a_i(\z)\abs{x_i-y_i}.
\EQN

\def\wk{\vk{k}}
\def\wl{\vk{l}}

Set for some $\delta\in(0,\frac{1}{2})$
\BQNY
&&\wk=(k_1\ldot   k_{n-1})\in \N^{n-1},\quad
 \wl=(l_1\ldot   l_{n-1})\in \N^{n-1},\quad J_{\wk}=\Pi_{j=1}^{n-1}[k_j\delta,(k_j+1)\delta],\\
&&\mathcal{L}_1=\LT\{\wk:\rho(J_{\wk}, L)\leq \frac{\delta}{3}\RT\},\quad
 \mathcal{L}_2=\LT\{\wk:J_{\wk}\subset L\RT\},\quad
D_{\wk}=\{\vk{x}: \abs{x_n-h(\wx)}\leq \delta,\ \wx\in J_{\wk}, \wk\in\mathcal{L}_1 \},\\
&&\mathcal{K}_1=\{(\wk,\wl):J_{\wk}\cap J_{\wl}\neq \emptyset, \wk,\wl\in\mathcal{L}_2 \},\quad
\mathcal{K}_2=\{(\wk,\wl):J_{\wk}\cap J_{\wl}= \emptyset,\wk,\wl\in\mathcal{L}_2\},\\
&& c_{\wk}=(k_1\delta\ldot  k_{n-1}\delta).
\EQNY
Here we need to notice that $\wk$ and $\wl$ are $(n-1)$-dimensional vector.\\
We have
\BQNY
\bigcup_{\wk\in \mathcal{L}_2 } J_{\wk} \subset L \subset\bigcup_{\wk\in \mathcal{L}_1 } J_{\wk}.
\EQNY
Bonferroni inequality leads to
\BQN
\pk{\sup_{\vk{x}\in E_0(\delta)}(X(\vk{x})+g(\vk{x}))>\mu}
&\leq& \sum_{\wk\in\mathcal{L}_1}\pk{\sup_{\vk{x}\in D_{\wk}}(X(\vk{x})+g(\vk{x}))>\mu},\label{Thm31upperbound1}\\
\pk{\sup_{\vk{x}\in E_0(\delta)}(X(\vk{x})+g(\vk{x}))>\mu}
&\geq& \sum_{\wk\in\mathcal{L}_2}\pk{\sup_{\vk{x}\in D_{\wk}}(X(\vk{x})+g(\vk{x}))>\mu}-\sum_{i=1}^2\Lambda_i(u),\label{Thm31lowerbound2}
\EQN
where
\BQNY
&&\Lambda_1(u)=\sum_{(\wk,\wl)\in\mathcal{K}_1}\pk{\sup_{\vk{x}\in D_{\wk}}(X(\vk{x})+g(\vk{x}))>\mu,\sup_{\vk{x}\in D_{\wl}}(X(\vk{x})+g(\vk{x}))>\mu},\\
&&\Lambda_2(u)=\sum_{(\wk,\wl)\in\mathcal{K}_2}\pk{\sup_{\vk{x}\in D_{\wk}}(X(\vk{x})+g(\vk{x}))>\mu,\sup_{\vk{x}\in D_{\wl}}(X(\vk{x})+g(\vk{x}))>\mu}.
\EQNY
By \eqref{Thm31sig}--\eqref{Thm31r3} and \nelem{lem1}, we have
\BQN\label{Thm31p3}
&&\sum_{\wk\in\mathcal{L}_1}\pk{\sup_{\vk{x}\in D_{\wk}}(X(\vk{x})+g(\vk{x}))>\mu}\nonumber\\
&&\sim\sum_{\wk\in\mathcal{L}_1}2^n\delta^{n-1} \LT(\Pi_{i=1}^{n-1}\tla_i(c_{\wk})\RT) \sqrt{\frac{\pi}{2}}e^{\frac{(2w)^2}{8}} \mu^{2n-1}\Psi(\mu)\nonumber\\
&&\sim\sum_{\wk\in\mathcal{L}_1}2^{3(n-1)}\delta^{n-1} \LT(\Pi_{i=1}^{n-1}\tla_i(c_{\wk})\RT) u^{2(n-1)}e^{-2u^2+2wu}\nonumber\\
&&\sim2^{3(n-1)} \int_{L} \LT(\Pi_{i=1}^{n-1}a_i(\wx, h(\wx))\RT)d \wx  u^{2(n-1)}e^{-2u^2+2uw},
\EQN
as $\mu\rw\IF, \delta\rw 0.$
\eE{Since $a_i(\wx,h(\wx)),\ \wx\in L$ is positive continuous and
	$L$ is Jordan measurable with positive Lebesgue measure, then $\int_{L} \LT(\Pi_{i=1}^{n-1}a_i(\wx, h(\wx))\RT)d \wx \in(0,\IF)$.}\\
Similarly,
\BQN\label{Thm31p4}
\sum_{\wk\in\mathcal{L}_2}\pk{\sup_{\vk{x}\in D_{\wk}}(X(\vk{x})+g(\vk{x}))>\mu}\sim2^{3(n-1)} \int_{L} \LT(\Pi_{i=1}^{n-1}a_i(\wx, h(\wx))\RT)d \wx  u^{2(n-1)}e^{-2u^2+2uw},
\EQN
as $\mu\rw\IF, \delta\rw 0.$
Next we will show that $\Lambda_i(u), i=1,2$ as $u\rw\IF$ are both negligible compared with
$$\sum_{\wk\in\mathcal{L}_2}\pk{\sup_{\vk{x}\in D_{\wk}}(X(\vk{x})+g(\vk{x}))>\mu}.$$
For any $(\wk,\wl)\in \mathcal{K}_1(u)$, without loss of generality, we assume that
$k_1+1=l_1$. Let
\BQNY
D_{\wk}^1=\LT[k_1\delta,(k_1+1)\delta-\delta^2\RT]
\times \Pi_{j=2}^{n-1}[k_j\delta,(k_j+1)\delta], \quad D_{\wk}^2=\LT[(k_1+1)\delta-\delta^2,(k_1+1)\delta\RT]
\times \Pi_{j=2}^{n-1}[k_j\delta,(k_j+1)\delta].
\EQNY
Then
\BQNY
&&\pk{\sup_{\vk{x}\in D_{\wk}}(X(\vk{x})+g(\vk{x}))>\mu,
\sup_{\vk{x}\in D_{\wl}}(X(\vk{x})+g(\vk{x}))>\mu}\\
&&\leq\pk{\sup_{\vk{x}\in D^1_{\wk}(u)}(X(\vk{x})+g(\vk{x}))>\mu,
\sup_{\vk{x}\in D_{\wl}(u)}(X(\vk{x})+g(\vk{x}))>\mu}+\pk{\sup_{\vk{x}\in D^2_{\wk}}(X(\vk{x})+g(\vk{x}))>\mu}.
\EQNY
Analogously as in \eqref{Thm31p3}, we have
\BQNY
\Lambda_{11}(u)&:=&\sum_{\wk\in\mathcal{L}_2}\pk{\sup_{\vk{x}\in D^2_{\wk}}(X(\vk{x})+g(\vk{x}))>\mu}\\
&\sim&\sum_{\wk\in\mathcal{L}_2}2^n\delta^{n} \LT(\Pi_{i=1}^{n-1}\tla_i(c_{\wk})\RT) \sqrt{\frac{\pi}{2}}e^{\frac{(2w)^2}{8}} \mu^{2n-1}\Psi(\mu)\nonumber\\
&\sim&\sum_{\wk\in\mathcal{L}_2}2^{3(n-1)}\delta^{n} \LT(\Pi_{i=1}^{n-1}\tla_i(c_{\wk})\RT)u^{2(n-1)}e^{-2u^2+2wu}\nonumber\\
&=&o\LT(u^{2(n-1)}e^{-2u^2+2uw}\RT),
\EQNY
as $u\rw\IF, \delta\rw 0.$\\
Moreover, since for $(\vk{x},\vk{y})\in D^1_{\wk}(u)\times D_{\wl}(u),\ (\wk,\wl)\in \mathcal{K}_1$, by \eqref{Thm31r2} we have
\BQNY
\Var(X(\vk{x})+X(\vk{y}))=2+2r(\vk{x},\vk{y})\leq 4-\mathbb{Q}\delta,
\EQNY
then by Borell-TIS inequality for $\mathbb{Q}=\sup_{\vk{x}\in [0,1]^n}g(\vk{x})$
\BQN\label{Thm31P5}
\Lambda_{12}(u)&:=&\sum_{(\wk,\wl)\in\mathcal{K}_1}\pk{\sup_{\vk{x}\in D^1_{\wk}(u)}(X(\vk{x})+g(\vk{x}))>\mu,
\sup_{\vk{x}\in D_{\wl}(u)}(X(\vk{x})+g(\vk{x}))>\mu}\nonumber\\
&\leq& \sum_{(\wk,\wl)\in\mathcal{K}_1}\pk{\sup_{(\vk{x},\vk{y})\in D^1_{\wk}(u)\times D_{\wl}(u) }(X(\vk{x})+X(\vk{y}) )>2(\mu-\mathbb{Q})}\nonumber\\
&\leq& \sum_{(\wk,\wl)\in\mathcal{K}_1}
\exp\LT(-\frac{(2(\mu-\mathbb{Q})-\mathbb{Q})^2}{2(4-\mathbb{Q}\delta)}\RT)\nonumber\\
&=&o\LT(u^{2(n-1)}e^{-2u^2+2uw}\RT),\ u\rw\IF.
\EQN
Since $J_{\wk}(u)$ has at most $3^{n-1}-1$ neighbors, then
\BQN\label{Thm31P6}
\Lambda_1(u)\leq 2(3^{n-1}-1)(\Lambda_{11}(u)+\Lambda_{12}(u))=o\LT(u^{2(n-1)}e^{-2u^2+2uw}\RT),\ u\rw\IF,\ \delta\rw0.
\EQN
Similarly, since for $(\vk{x},\vk{y})\in D_{\wk}(u)\times D_{\wl}(u),\ (\wk,\wl)\in \mathcal{K}_2$, by \eqref{Thm31r2} we have
\BQNY
\Var(X(\vk{x})+X(\vk{y}))=2+2r(\vk{x},\vk{y})\leq 4-\mathbb{Q}\delta,
\EQNY
and then we have
\BQN\label{Thm31P7}
\Lambda_2(u)&=&\sum_{(\wk,\wl)\in\mathcal{K}_2}\pk{\sup_{\vk{x}\in D_{\wk}}(X(\vk{x})+g(\vk{x}))>\mu,\sup_{\vk{x}\in D_{\wl}}(X(\vk{x})+g(\vk{x}))>\mu}\nonumber\\
&\leq& \sum_{(\wk,\wl)\in\mathcal{K}_2}\pk{\sup_{(\vk{x},\vk{y})\in D_{\wk}(u)\times D_{\wl}(u) }(X(\vk{x})+X(\vk{y}) )>2(\mu-\mathbb{Q})}\nonumber\\
&\leq& \sum_{(\wk,\wl)\in\mathcal{K}_2}
\exp\LT(-\frac{(2(\mu-\mathbb{Q})-\mathbb{Q})^2}{2(4-\mathbb{Q}\delta)}\RT)\nonumber\\
&=&o\LT(u^{2(n-1)}e^{-2u^2+2uw}\RT), \ u\rw\IF.
\EQN

Inserting \eqref{Thm31p3}, \eqref{Thm31p4}, \eqref{Thm31P6}, and \eqref{Thm31P7} into \eqref{Thm31upperbound1} and \eqref{Thm31lowerbound2} implies
\BQNY
\pk{\sup_{\vk{x}\in E_0(\delta)}(X(\vk{x})+g(\vk{x}))>\mu}\sim 2^{3(n-1)} \int_{L} \LT(\Pi_{i=1}^{n-1}a_i(\wx, h(\wx))\RT)d \wx  u^{2(n-1)}e^{-2u^2+2uw},\ u\rw\IF,
\EQNY
which combined with \eqref{Thm31ub1} and
\eqref{Thm31lb2} establishes the claim  \eqref{Mulre1}.

\QED

\proofkorr{Corr3}
Since we assume that $f$ is positive in $[0,1]^n$, we have for any $\x \in (0,1)^{n-1}$ the $i$th partial derivative of $F$ denoted by $a_i(\x)$ is positive
and continuous. Let $\Q=\sup_{\x\in [0,1]^n} f(\x)$ which is finite and positive by the assumption and $\seE=\Bigl\{\vk{x}\in [0,1]^n:F(\vk{x})=\frac{1}{2}\Bigr\}$. Using Taylor expansion, we have
\BQNY
\Abs{F(\x)-F(\z)-\sum_{i=1}^{n}a_i(\vk{z})(x_i-z_i)}
\leq \Q\LT(\sum_{i=1}^{n}\abs{x_i-z_i}^2\RT)
\EQNY
for any $\x, \z \in [0,1]^{n}$. Consequently,
\BQNY
\sup_{\z\in \seE}\sup_{0 < \abs{\x-\z}<\delta}\frac{\abs{F(\x)-F(\z)-\sum_{i=1}^na_i(\vk{z})(x_i-z_i)}}
{\sum_{i=1}^n\abs{x_i-z_i}}\leq 2\Q\delta,
\EQNY
which combined with the continuity of $a_i$  implies
\BQNY
\lim_{\delta\rw 0}\sup_{\vk{z}\in \seE} \underset{\vk{x}\neq \vk{y}}{\sup_{\abs{\vk{x}-\vk{z}},\abs{\vk{y}-\vk{z}}\leq \delta}}\frac{\abs{ F(\vk{x})-F(\vk{y})-\sum_{i=1}^{n}a_i(\vk{z})(x_i-y_i)}}
{\sum_{i=1}^{n}\abs{x_i-y_i}}=0.
\EQNY
In view of \netheo{ThDK} the set $\seE$ is not empty and moreover its projection on  $[0,1]^{n-1}$  denoted by  $L$  is Jordan measurable with positive Lebesgue measure (with respect to $\lambda_{n-1})$.\\
Set $L^o$ is the interior of $L$. By the positivity of the partial derivatives on interior of $[0,1]^n$ and the fact that $F$ is strictly increasing on $[0,1]^{n-1}$, we have that for any $\wx \in L^o$ there is only one $x_n$ such that $F(\wx,x_n)=1/2$. Consequently, $x_n= g(\wx)$ for some bijective function $g$ for any $\wx\in L^o$. Since $F$ is continuously differentiable on $L^o$, for any $\x \in \seE$ with $\wx\in L^o$ by the implicit function theorem there exists $\ve>0$ such that for any $\wy \in \mathcal{O}_\vn(\wx)=\{\wz\in I^{n-1}:\abs{\wz-\wx}<\vn\}\subseteq L^o$, we have
$ F(\wy, h_{\x}(\wy))=1/2$. By the above, $\lambda_{n-1}(L^o)\geq \lambda_{n-1}(\mathcal{O}_\vn(\wx) )>0$ and $h_{\x}$ does not depend on $\x$ and $h_{\x}(\wx)=g(\wx)$ for any $\wx\in L^o$. It follows thus that $g$ is continuously differentiable on $L^o$. Moreover for any $\wx\in L^o$
\BQNY
\frac{\partial h(\wx)}{\partial x_i }=-\frac{a_i(\wx,h(\wx))}{a_n(\wx,h(\wx))}<0,\ i=1\ldot  n-1,
\EQNY
and thus $h$ is continuously differentiable in $L^o$. Hence the proof follows by \netheo{ThmM1} since \eqref{FEXM}, the  continuous differentiability of $h$ in $L^o$ and the Jordan measurability of $L$ are satisfied.

\QED

\prooftheo{Thm2} The variance function of $W_F(\vk{x})$ is
\BQNY
\sigma^2_F(\vk{x})=\Var(W_F(\vk{x}))=F(\vk{x})(1-F(\vk{x})),\quad \vk{x}\in [0,1]^n,
\EQNY
which attains its maximum equal to $\frac{1}{4}$ over $[0,1]^n$ at $(\vk{z})$ which satisfies $F(\vk{z})=\frac{1}{2}$.\\
Since $\sigma^2_F(\vk{x}), \vk x \in [0,1]^n$ is a continuous function, we have
\BQNY
\sup_{\vk{x}\in [0,1]^n\setminus \seEd}\sigma^2_F(\vk{x})=\frac{1}{4}-\delta^2.
\EQNY
By Borell-TIS inequality (ref.\cite{AdlerTaylor}), as $u\rw\IF$
\BQN\label{boundth21}
\pk{\sup_{\vk{x}\in [0,1]^n\setminus \seEd}W(\vk{x})>u\Big|W(\vk{1})=w}&=&\pk{\sup_{\vk{x}\in [0,1]^n\setminus \seEd}\LT(W_F(\vk{x})+F(\vk{x})w\RT)>u}\nonumber\\
&\leq&\pk{\sup_{\vk{x}\in [0,1]^n\setminus \seEd}W_F(\vk{x})>u-\mathbb{Q}_1}\nonumber\\
&\leq&\exp\LT(-\frac{(u-\mathbb{Q}_1-\mathbb{Q}_2)^2}{2\sigma^2_m}\RT)\nonumber\\
&=&o\LT(e^{-2u+2uw}\RT),
\EQN
where $\mathbb{Q}_1=\sup_{\vk{x}\in [0,1]^n}F(\vk{x})$ and
$\mathbb{Q}_2=\E{\sup_{\vk{x}\in [0,1]^n}W_F(\vk{x})}<\IF$.
By  \eqref{condition1}, we know that $W_F(\vk{x}), \vk{x}\in \seEd$ is a Gaussian fields with covariance function
\BQN\label{covariance1}
\E{W_F(\vk{x})W_F(\vk{y})}&=& F(\vk{x}\wedge \vk{y})-F(\vk{x})F(\vk{y})\nonumber\\
&=&F(\vk{x})\wedge F(\vk{y})-F(\vk{x})F(\vk{y})\nonumber\\
&=&\E{\eE{B_0}(F(\vk{x}))B_0(F(\vk{y}))}
\EQN
where $B_0(t)=B(t)-tB(1)$ is the standard Brownian bridge.
Then by \cite{GauTrend16} [Example 3.12]
\BQNY
\pk{\sup_{\vk{x}\in\seEd}W(\vk{x})>u\Big|W(\vk{1})=w}&=&\pk{\sup_{\vk{x}\in \seEd}\LT(W_F(\vk{x})+F(\vk{x})w\RT)>u}\\
&=&\pk{\sup_{\vk{x}\in \seEd}\LT(B_0(F(\vk{x}))+F(\vk{x})w\RT)>u}\\
&=&\pk{\sup_{F(\vk{x})\in [\frac{1}{2}-\delta,\frac{1}{2}+\delta]}\LT(B_0(F(\vk{x}))+F(\vk{x})w\RT)>u}\\
&=&\pk{\sup_{x\in [\frac{1}{2}-\delta,\frac{1}{2}+\delta]}\LT(\eE{B_0}(x)+wx\RT)>u}\\
&\sim&  e^{-2u^2+2uw}, \ u\rw\IF,
\EQNY
which combined with \eqref{boundth21} and the fact that
\begin{align*}
\pk{\sup_{\vk{x}\in [0,1]^n}W(\vk{x})>u\Big|W(\vk{1})=w}&\geq \pk{\sup_{\vk{x}\in\seEd}W(\vk{x})>u\Big|W(\vk{1})=w}\\
 \pk{\sup_{\vk{x}\in [0,1]^n}W(\vk{x})>u\Big|W(\vk{1})=w}&\leq \pk{\sup_{\vk{x}\in\seEd}W(\vk{x})>u\Big|W(\vk{1})=w}
+\pk{\sup_{\vk{x}\in [0,1]^n\setminus \seEd}W(\vk{x})>u\Big|W(\vk{1})=w}
\end{align*}
implies
\BQN\label{re1}
\pk{\sup_{\vk{x}\in [0,1]^n}W(\vk{x})>u\Big|W(\vk{1})=w}\sim e^{-2u^2+2uw}, \ u\rw\IF.
\EQN

If \eqref{condition1}  holds for $\vk{x}, \vk{y}\in [0,1]^n$, then by \eqref{covariance1} for any $u>0$
\BQNY
\pk{\sup_{\vk{x}\in [0,1]^n}W(\vk{x})>u\Big|W(\vk{1})=w}&=&\pk{\sup_{\vk{x}\in [0,1]^n}\LT(W_F(\vk{x})+F(\vk{x})w\RT)>u}\\
&=&\pk{\sup_{F(\vk{x})\in [0,1]}\LT(B_0(F(\vk{x}))+F(\vk{x})w\RT)>u}\\
&=&\pk{\sup_{x\in [0,1]}\LT(B_0(x)+wx\RT)>u}\\
&=&  e^{-2u^2+2uw},
\EQNY
where the last equation is well-known,  see e.g., \cite{AsymBB2003} [Lemma 2.7].

\QED

\prooftheo{Corr4} For $u>0$ we have
\begin{align*}
\pk{\sup_{\vk{x}\in \mathcal{E}}\abs{W(\vk{x})}>u\Big|W(\vk{1})=w}
&=\pk{\sup_{\vk{x}\in  \mathcal{E}}W(\vk{x})>u\Big|W(\vk{1})=w}
+\pk{\inf_{\vk{x}\in  \mathcal{E}}W(\vk{x})<-u\Big|W(\vk{1})=w}\nonumber\\
&\quad-\pk{\sup_{\vk{x}\in  \mathcal{E}}W(\vk{x})>u\Big|W(\vk{1})=w
,\ \inf_{\vk{x}\in  \mathcal{E}}W(\vk{x})<-u\Big|W(\vk{1})=w}\nonumber\\
&=: J_1(u)+J_2(u)-J_3(u).
\end{align*}
Since there exist $\x_0\in \mathcal{E}$ such that $F(\x_0)=\frac{1}{2}$ and
\BQNY
\Var(W(\x_0)-F(\x_0)W(1))=F(\x_0)-F^2(\x_0)=\frac{1}{4},
\EQNY
we have
\BQNY
J_1(u)=\pk{\sup_{\vk{x}\in  \mathcal{E}}\LT(W(\vk{x})-F(\x)W(\vk{1})+F(\x) w\RT)>u}
\geq \pk{W(\x_0)-F(\x_0)W(\vk{1})+\frac{1}{2} w>u}=\Psi(2u-w),
\EQNY
and
\BQNY
J_2(u)=\pk{\inf_{\vk{x}\in  \mathcal{E}}\LT(W(\vk{x})-F(\x)W(\vk{1})+F(\x) w\RT)<-u}
\geq \pk{W(\x_0)-F(\x_0)W(\vk{1})+\frac{1}{2} w<-u}=\Psi(2u+w).
\EQNY
Thus we have for $u>0$
\BQNY
J_1(u)+J_2(u)\geq \Psi(2u-w)+\Psi(2u+w)\geq \Psi(2u) .
\EQNY
Next in order to get the finial result, we need to show that
\BQNY
J_3(u)=o\LT(\Psi(2u)\RT), \ u\rw\IF.
\EQNY
We have that for $u>0$
\BQNY
J_3(u)&\leq&\pk{\sup_{\vk{x}\in \mathcal{E}\setminus \seEd}W(\vk{x})>u\Big|W(\vk{1})=w
}+\pk{\inf_{\y\in \mathcal{E}\setminus \seEd }W(\vk{y})<-u\Big|W(\vk{1})=w}\\
&&+\pk{\sup_{\vk{x}\in \seEd }W(\vk{x})>u\Big|W(\vk{1})=w
,\ \inf_{\y\in \seEd}W(\y)<-u\Big|W(\vk{1})=w}\\
&=:& J_{31}(u)+J_{32}(u)+J_{33}(u).
\EQNY
Since for $\x\in \mathcal{E}\setminus \seEd $, $\abs{F(\x)-\frac{1}{2}}>\delta$ and
\BQNY
\sigma^2_m:=\sup_{\vk{x}\in \mathcal{E}\setminus \seEd}\Var\LT(W(\vk{x})-F(\vk{x})W(\vk{1})\RT)=\sup_{\vk{x}\in \mathcal{E}\setminus \seEd}F(\x)(1-F(\x))< \frac{1}{4}-\delta^2.
\EQNY
By Borell-TIS inequality (ref.\cite{AdlerTaylor}), we have for
all $u$ sufficiently large
\BQNY
J_{31}(u)&=&\pk{\sup_{\vk{x}\in \mathcal{E}\setminus \seEd}\LT(W(\vk{x})-F(\vk{x})W(\vk{1})+F(\vk{x})w\RT)>u}\\
&\leq&\pk{\sup_{\vk{x}\in \mathcal{E}\setminus \seEd}\LT(W(\vk{x})-F(\vk{x})W(\vk{1})\RT)>u-\abs{w}}\\
&\leq&\exp\LT(-\frac{(u-\abs{w}-\mathbb{Q}_1)^2}{2\sigma^2_m}\RT)
=o\LT(\Psi(2u)\RT),
\EQNY
and
\BQNY
J_{32}(u)&=&\pk{\inf_{\y\in \mathcal{E}\setminus \seEd }\LT(W(\vk{x})-F(\vk{x})W(\vk{1})+F(\vk{x})w\RT)<-u}\\
&=&\pk{\sup_{\vk{x}\in \mathcal{E}\setminus \seEd }\LT(-W(\vk{x})+F(\vk{x})W(\vk{1})-F(\vk{x})w\RT)>u}\\
&\leq&\pk{\sup_{\vk{x}\in \mathcal{E}\setminus \seEd }\LT(W(\vk{x})-F(\vk{x})W(\vk{1})\RT)>u-\abs{w}}\\
&\leq& \exp\LT(-\frac{(u-\abs{w}-\mathbb{Q}_1)^2}{2\sigma^2_m}\RT)
=o\LT(\Psi(2u)\RT),
\EQNY
where we use the symmetry of $(W(\vk{x})-F(\vk{x})W(\vk{1}))$ and $\mathbb{Q}_1:=\sup_{\x\in [0,1]^n\setminus \seEd}\E{W(\vk{x})-F(\vk{x})W(\vk{1})}\in(0,\IF)$.
Further,  by \eqref{Con2}
\begin{align*}
\varrho:&=\sup_{(\vk{x},\vk{y})\in \seEd\times \seEd}\Var\LT((W(\vk{x})-W(\vk{y}))-(F(\vk{x})-F(\vk{y}))W(\vk{1})\RT)\\
&=\sup_{(\vk{x},\vk{y})\in \seEd\times \seEd}\LT(F(\vk{x})+F(\vk{y})-2F(\vk{x}\wedge\vk{y})-(F(\vk{x})-F(\vk{y}))^2\RT)\\
&\leq \sup_{(\vk{x},\vk{y})\in \seEd\times \seEd}\LT(F(\x)+F(\y)-F^2(\x)-F^2(\y)+2F(\vk{x})F(\vk{y})\RT)-
\inf_{(\vk{x},\vk{y})\in \seEd\times \seEd}2F(\vk{x}\wedge\vk{y})\\
&< 1+2\delta-2\delta=1,
\end{align*}
where we use the fact that
$$\sup_{(x,y)\in[\frac{1}{2}-\delta,\frac{1}{2}+\delta]
\times[\frac{1}{2}-\delta,\frac{1}{2}+\delta]}\LT(x+y-x^2-y^2+2xy\RT)=1+2\delta.$$
By Borell-TIS inequality (ref.\cite{AdlerTaylor}) again
\BQN\label{Thm31pp}
J_{33}(u)&=&\pk{\sup_{\vk{x}\in  \seEd }\LT(W(\x)-F(\x)W(\vk{1})+F(\x) w\RT)>u
,\ \inf_{\y\in  \seEd}\LT(W(\y)-F(\y)W(\vk{1})+F(\y) w\RT)<-u}\nonumber\\
&=&\pk{\sup_{\vk{x}\in  \seEd }\LT(W(\x)-F(\x)W(\vk{1})+F(\x) w\RT)>u
,\ \sup_{\y\in  \seEd}\LT(-W(\y)+F(\y)W(\vk{1})-F(\y) w\RT)>u}\nonumber\\
&\leq&\pk{\sup_{(\vk{x},\vk{y})\in  \seEd\times \seEd}\LT((W(\vk{x})-W(\vk{y}))-(F(\vk{x})-F(\vk{y}))W(\vk{1})
+(F(\vk{x})-F(\vk{y}))w\RT)>2u}\nonumber\\
&\leq& \pk{\sup_{(\vk{x},\vk{y})\in  \seEd\times \seEd}\LT((W(\vk{x})-W(\vk{y}))-(F(\vk{x})-F(\vk{y}))W(\vk{1})
\RT)>2u-\abs{w}}\nonumber\\
&\leq& \exp\LT(-\frac{\LT(2u-\abs{w}-\mathbb{Q}_2\RT)^2}{2\varrho}\RT)
= o\LT(\Psi(2u)\RT),
\EQN
where $\mathbb{Q}_2=\sup_{(\vk{x},\vk{y})\in \seEd\times \seEd}\LT((W(\vk{x})-W(\vk{y}))-(F(\vk{x})-F(\vk{y}))W(\vk{1})\RT)<\IF$.
Hence, as $u\rw\IF$
\BQN\label{JJ3}
J_3(u)\leq J_{31}(u)+J_{32}(u)+J_{33}(u)=o\LT(\Psi(2u)\RT),
\EQN
and the proof is complete.

\QED

\proofprop{PROP00}
i) We have that
\BQNY
\seE=\LT\{\vk{x}\in [0,1]^2: x_1+x_2=\frac{3}{2}\RT\}, \quad L=\LT[\frac{1}{2},1\RT]
\EQNY
and further, for $\vk{x},\vk{y}\in \LT\{\z\in[0,1]^2:z_1+z_2\geq 1\RT\}\supset\seE$
\BQNY
F(\vk{x})-F(\vk{y})=\sum_{i=1}^2(x_i-y_i)
\EQNY
which implies \eqref{FEXM} are satisfied.\\
In view of \netheo{ThmM1}, we have $K=4$ by taking  $a_1(x,h(x))=1, \ x \in L=\LT[\frac{1}{2},1\RT]$.

For $\delta\in\LT(0,\frac{1}{4}\RT)$ small enough, set $\mathcal{E}=\LT[\frac{1}{2}+\delta, 1\RT]\times\LT[\frac{1}{2}+\delta, 1\RT]$ and $\seEd=\LT\{\x\in \mathcal{E}: \frac{1}{2}-\delta \leq F(\vk{x})\leq \frac{1}{2}+\delta \RT\}$, then we have
\BQNY
 \inf_{\x,\y\in\seEd} F(\x\wedge\y)= \inf_{\x,\y\in\seEd} (x_1\wedge y_1+x_2\wedge y_2-1)
 \geq 2\delta>\delta,
\EQNY
which show that \eqref{Con2} holds.
We have for $u>0$
\BQN
&&\pk{\sup_{\vk{x}\in [0,1]^2}\abs{W(\vk{x})}>u\Big|W(\vk{1})=w}\geq J_1(u),\label{lowbound1}\\
&&\pk{\sup_{\vk{x}\in [0,1]^2}\abs{W(\vk{x})}>u\Big|W(\vk{1})=w}\leq\sum_{i=1}^3 J_i(u),\label{upbound1}
\EQN
where
\BQNY
&&J_1(u)=\pk{\sup_{\vk{x}\in \mathcal{E}}\abs{W(\vk{x})}>u\Big|W(\vk{1})=w},
\ \ J_2(u)=\pk{\sup_{\vk{x}\in \LT[0,\frac{1}{2}+\delta\RT]\times [0,1]}\abs{W(\vk{x})}>u\Big|W(\vk{1})=w}\\
&&J_3(u)=\pk{\sup_{\vk{x}\in[0,1]\times\LT[0,\frac{1}{2}+\delta\RT]}
\abs{W(\vk{x})}>u\Big|W(\vk{1})=w}.
\EQNY

For $J_1(u)$, by \netheo{ThmM1} with $a_1(x,h(x))=1, \ x \in L_1=\LT[\frac{1}{2}+\delta,1\RT]$ and \netheo{Corr4}
\BQNY
J_1(u)&\sim& \pk{\sup_{\vk{x}\in \mathcal{E}}W(\vk{x})>u\Big|W(\vk{1})=w}+\pk{\sup_{\vk{x}\in \mathcal{E}}W(\vk{x})>u\Big|W(\vk{1})=-w} \\
&\sim& c(4-8\delta)u^{2}e^{-2u^2+2u\abs{w}} \\
&\sim& 4cu^{2}e^{-2u^2+2u\abs{w}},\ u\rw\IF,\ \delta\rw 0
\EQNY
holds with $c=1$ for $w\neq 0$ and $c=2$ for $w=0$.
Since  by \netheo{ThmM1} with $a_1(x,h(x))=1, \ x \in L_2=\LT[\frac{1}{2},\frac{1}{2}+\delta\RT]$
\BQNY
J_2(u)&\leq&\pk{\sup_{\vk{x}\in \LT[0,\frac{1}{2}+\delta\RT]\times [0,1]}W(\vk{x})>u\Big|W(\vk{1})=w}+\pk{\inf_{\vk{x}\in \LT[0,\frac{1}{2}+\delta\RT]\times [0,1]}W(\vk{x})<-u\Big|W(\vk{1})=w}\\
&\leq&\pk{\sup_{\vk{x}\in \LT[0,\frac{1}{2}+\delta\RT]\times [0,1]}W(\vk{x})>u\Big|W(\vk{1})=w}+\pk{\sup_{\vk{x}\in \LT[0,\frac{1}{2}+\delta\RT]\times [0,1]}W(\vk{x})>u\Big|W(\vk{1})=-w}\\
&\sim& 8\delta u^{2}e^{-2u^2+2uw}+8\delta u^{2}e^{-2u^2-2uw}\\
&\sim& o\LT(u^{2}e^{-2u^2+2u\abs{w}}\RT),\ u\rw\IF,\ \delta\rw 0.
\EQNY
Similarly,
\BQNY
J_3(u)=o\LT(u^{2}e^{-2u^2+2u\abs{w}}\RT),\ u\rw\IF,\ \delta\rw 0.
\EQNY
Thus  by \eqref{lowbound1} and \eqref{upbound1}
\BQNY
\pk{\sup_{\vk{x}\in [0,1]^2}\abs{W(\vk{x})}>u\Big|W(\vk{1})=w}\sim J_1(u)\sim 4cu^{2}e^{-2u^2+2u\abs{w}},\ u\rw\IF.
\EQNY
ii) We have that
\BQNY
\seE=\LT\{\vk{x}\in [0,1]^2: x_1+x_2+x_1x_2=2\RT\}, \quad L=\LT[\frac{1}{2},1\RT],
\ h(x_1)=\frac{2-x_1}{1+x_1},
\EQNY
and for a small neighborhood of $\seE$, the density function of $F(x_1,x_2)$
\BQNY
f(x_1,x_2)= \frac{4-2x_1-2x_2}{(1+(1-x_1)(1-x_2))^3}>0.
\EQNY
In view of \nekorr{Corr3}, we have $K=3\ln 3$ by taking
$a_1(x,h(x))=\frac{3}{4-4(1-x)^2}, \ x \in L=\LT[\frac{1}{2},1\RT]$.

For $\delta\in(0,\frac{3-2\sqrt{2}}{2})$ , set $\seEd=\LT\{\x\in [0,1]^2: \frac{1}{2}-\delta \leq F(\vk{x})\leq \frac{1}{2}+\delta \RT\}$, then we have
\BQNY
\min(x_1,x_2)\geq \frac{1}{2}-\delta,\ \x\in\seEd,
\EQNY
and
\BQNY
 \inf_{\x,\y\in\seEd} F(\x\wedge\y)= \inf_{\x,\y\in\seEd} \frac{(x_1\wedge y_1)\times(x_2\wedge y_2-1)}{1+(1-x_1\wedge y_1)(1-x_2\wedge y_2)}
 \geq \frac{(\frac{1}{2}-\delta)^2}{1+(\frac{1}{2}+\delta)^2}>
 \frac{(\frac{1}{2}-\delta)^2}{2}>\delta,
\EQNY
which show that \eqref{Con2} holds.
Thus by \netheo{ThmM1} and \netheo{Corr4}
\BQNY
\pk{\sup_{\vk{x}\in [0,1]^2}\abs{W(\vk{x})}>u\Big|W(\vk{1})=w}
\sim c K u^{2}e^{-2u^2+2u\abs{w}}
\EQNY
holds with $c=1$ for $w\neq 0$ and $c=2$ for $w=0$

\QED

\proofprop{PROP0}
i) We have that
\BQNY
\seE=\LT\{\vk{x}\in [0,1]^n:\Pi_{i=1}^n x_i=\frac{1}{2}\RT\}, \quad L=\LT\{\wx\in [0,1]^{n-1}:\frac{1}{2}\leq \Pi_{i=1}^{n-1}x_i\leq 1\RT\}.
\EQNY
Further, $F(\x)$ has the density function $f(\x)\equiv1$.
By \nekorr{Corr3}, we get the results with
$$a_i(\wx,h(\wx))=\frac{1}{2x_i}, \quad  \wx \in L,\ i=1\ldot n-1.$$
Further, if we take $\mathcal{E}=[0,1]^n$ and $\seEd:=\LT\{\x\in [0,1]^n: \frac{1}{2}-\delta \leq F(\vk{x})\leq \frac{1}{2}+\delta \RT\}$ with $\delta\in(0,\frac{1}{4})$ such that $\delta^{1/n}+\delta<\frac{1}{2}$,
then we have
$$\min_{1\leq i\leq n} x_i\geq \frac{1}{2}-\delta , \ \x \in \seEd,$$
and
\BQNY
 \inf_{\x,\y\in\seEd} F(\x\wedge\y)= \inf_{\x,\y\in\seEd} \Pi_{i=1}^n(x_i\wedge y_i)
 \geq \LT(\frac{1}{2}-\delta\RT)^n>\delta,
\EQNY
which show that \eqref{Con2} holds. Thus by \nekorr{Corr3} and \netheo{Corr4}
\BQNY
&&\pk{\sup_{\vk{x}\in [0,1]^n}\abs{W(\vk{x})}>u\Big|W(\vk{1})=w}\nonumber\\
&&\sim \pk{\sup_{\vk{x}\in [0,1]^n}W(\vk{x})>u\Big|W(\vk{1})=w}+\pk{\sup_{\vk{x}\in [0,1]^n}W(\vk{x})>u\Big|W(\vk{1})=-w}\\
&&\sim c K u^{2(n-1)}e^{-2u^2+2u\abs{w}}
\EQNY
holds with $c=1$ for $w\neq 0$ and $c=2$ for $w=0$.

\COM{ii) We have that
\BQNY
\seE=\LT\{\vk{x}\in [0,1]^n:\sum_{i=1}^n x_i=\frac{n}{2}\RT\}, \quad L=\LT\{\wx\in [0,1]^{n-1}:\frac{n-2}{2}\leq \sum_{i=1}^{n-1}x_i\leq \frac{n}{2}\RT\}.
\EQNY
Further, for $\vk{x},\vk{y}\in [0,1]^n$
\BQNY
F(\vk{x})-F(\vk{y})=\frac{1}{n}\sum_{i=1}^n(x_i-y_i).
\EQNY
In the notation of  \netheo{ThmM1}, we have  $a_i(\wx,h(\wx))=\frac{1}{n}, \ \wx \in L$, hence the claim follows from the aforementioned theorem. \\}

ii) For $\delta>0$ small enough, set
\BQNY
&&E_i(\delta)=\{\x\in [0,1]^n: x_i\leq \min(x_1\ldot x_{i-1},x_{i+1}\ldot x_n)-\delta\},\
i=1\ldot n,\\
&&E_{n+1}(\delta)=[0,1]^n\setminus \LT(\cup_{i=1}^n E_i(\delta)\RT).
\EQNY
We have
\begin{align*}
\pk{\sup_{\vk{x}\in [0,1]^n}\LT(W(\vk{x})\Big|W(\vk{1})=w\RT)>u}&\leq
\sum_{i=1}^{n+1}\pk{\sup_{\vk{x}\in E_i(\delta)}\LT(W(\vk{x})\Big|W(\vk{1})=w\RT)>u},\\
\pk{\sup_{\vk{x}\in [0,1]^n}\LT(W(\vk{x})\Big|W(\vk{1})=w\RT)>u}&\geq
\sum_{i=1}^{n}\pk{\sup_{\vk{x}\in E_i(\delta)}\LT(W(\vk{x})\Big|W(\vk{1})=w\RT)>u}\\
&-\sum_{1\leq i<j\leq n}\pk{\sup_{\vk{x}\in E_i(\delta)}\LT(W(\vk{x})\Big|W(\vk{1})=w\RT)>u,
\sup_{\vk{x}\in E_j(\delta)}\LT(W(\vk{x})\Big|W(\vk{1})=w\RT)>u}.
\end{align*}
For $\x\in E_n(\delta)$, we have $F(\x)=x_n\LT(d+(1-d)\Pi_{i=1}^{n-1}x_i\RT)$,
\BQNY
&&\seE_n=\{\x\in E_n(\delta):F(\x)=\frac{1}{2}\}=\LT\{(\wx,h(\wx)):\wx\in
L\RT\},\\
&&L(\delta)=\LT\{\wx\in[0,1]^{n-1}:h(\wx)\leq\min_{1\leq i\leq n-1}x_i-\delta\RT\},\ h(\wx)=\frac{1}{2d+2(1-d)\Pi_{i=1}^{n-1}x_i},
\EQNY
and for $\x\in E_n(\delta)$, $F(\x)$ has density function $f(\x)\equiv1$.

Thus by \nekorr{Corr3}, we have as $u\rw\IF$
\BQNY
\pk{\sup_{\vk{x}\in E_n(\delta)}\LT(W(\vk{x})\Big|W(\vk{1})=w\RT)>u}
\sim (4d)^{(n-1)} \int_{L(\delta)} \frac{\LT(\Pi_{i=1}^{n-1}x_i\RT)^{n-2}}{\LT(d+(1-d)\Pi_{i=1}^{n-1}x_i\RT)^{n-1}}d \wx u^{2(n-1)}e^{-2u^2+2uw}.
\EQNY
Further, we have
\BQNY
\lim_{\delta\rw 0}\lim_{u\rw\IF}\frac{\pk{\sup_{\vk{x}\in E_n(\delta)}\LT(W(\vk{x})\Big|W(\vk{1})=w\RT)>u}}{u^{2(n-1)}e^{-2u^2+2uw}}
= (4d)^{(n-1)} \int_{L(\delta)} \frac{\LT(\Pi_{i=1}^{n-1}x_i\RT)^{n-2}}{\LT(d+(1-d)\Pi_{i=1}^{n-1}x_i\RT)^{n-1}}d \wx  ,
\EQNY
where
\BQNY
L=\LT\{\wx\in[0,1]^{n-1}:\frac{1}{2d+2(1-d)\Pi_{i=1}^{n-1}x_i}\leq\min_{1\leq i\leq n-1}x_i\RT\}.
\EQNY
By the symmetry of $F(\x)$ on $E_i(\delta), i=1\ldot n-1$,
\BQNY
\lim_{\delta\rw 0}\lim_{u\rw\IF}\frac{\pk{\sup_{\vk{x}\in E_{i}(\delta)}\LT(W(\vk{x})\Big|W(\vk{1})=w\RT)>u}}{u^{2(n-1)}e^{-2u^2+2uw}}
= (4d)^{(n-1)} \int_{L} \frac{\LT(\Pi_{i=1}^{n-1}x_i\RT)^{n-2}}{\LT(d+(1-d)\Pi_{i=1}^{n-1}x_i\RT)^{n-1}}d \wx  .
\EQNY
Assume that  $W_i(\x)$ is a Brownian sheet based on $F_i(\x)=d x_i +(1-d)\Pi_{j=1}^{n}x_j$, $i=1\ldot n$. We have
\BQNY
\pk{\sup_{\vk{x}\in E_{n+1}(\delta)}\LT(W(\vk{x})\Big|W(\vk{1})=w\RT)>u}
\leq \sum_{i=1}^n\pk{\sup_{\vk{x}\in E_{n+1}(\delta)}\LT(W_i(\vk{x})\Big|W_i(\vk{1})=w\RT)>u}.
\EQNY
For $F_n(\x)=d x_n+(1-d)\Pi_{i=1}^{n}x_i$ and $\x\in E_{n+1}(\delta)$ we have
\BQNY
&&\seE_{n+1}=\LT\{\x\in E_{n+1}(\delta):F_n(\x)=\frac{1}{2}\RT\}=\LT\{(\wx,h(\wx)):\wx\in
L_{n+1}\RT\},\\
&&L_{n+1}(\delta)=\LT\{\wx\in[0,1]^{n-1}:(\wx,h(\wx))\in\seE_{n+1}\RT\},\ h(\wx)=\frac{1}{2d+2(1-d)\Pi_{i=1}^{n-1}x_i}.
\EQNY
Further by \nekorr{Corr3}, we have
\BQNY
\pk{\sup_{\vk{x}\in E_{n+1}(\delta)}\LT(W_n(\vk{x})\Big|W_n(\vk{1})=w\RT)>u}
\sim  (4d)^{(n-1)}\int_{L_{n+1}(\delta)} \frac{\LT(\Pi_{i=1}^{n-1}x_i\RT)^{n-2}}{\LT(d+(1-d)\Pi_{i=1}^{n-1}x_i\RT)^{n-1}}d \wx u^{2(n-1)}e^{-2u^2+2uw}
, u\rw\IF.
\EQNY
Similarly, for $i=1\ldot n-1$
\BQNY
\pk{\sup_{\vk{x}\in E_{n+1}(\delta)}\LT(W_i(\vk{x})\Big|W_i(\vk{1})=w\RT)>u}
\sim  (4d)^{(n-1)}\int_{L_{n+1}(\delta)} \frac{\LT(\Pi_{i=1}^{n-1}x_i\RT)^{n-2}}{\LT(d+(1-d)\Pi_{i=1}^{n-1}x_i\RT)^{n-1}}d \wx u^{2(n-1)}e^{-2u^2+2uw}
, u\rw\IF.
\EQNY
Since $\lim_{\delta\rw 0}\LT(\cup_{i=1}^{n}E_i(\delta)\RT)=[0,1]^{n}$ and $\lim_{\delta\rw 0} E_{n+1}(\delta)=\emptyset$, then  $\lim_{\delta\rw 0}\lambda_{n-1}(L_{n+1}(\delta))=0$.
Thus we have
\BQNY
\pk{\sup_{\vk{x}\in E_{n+1}(\delta)}\LT(W(\vk{x})\Big|W(\vk{1})=w\RT)>u}
&\leq& n(4d)^{(n-1)}\int_{L_{n+1}(\delta)} \frac{\LT(\Pi_{i=1}^{n-1}x_i\RT)^{n-2}}{\LT(d+(1-d)\Pi_{i=1}^{n-1}x_i\RT)^{n-1}}d \wx u^{2(n-1)}e^{-2u^2+2uw}\\
&=&o\LT(u^{2(n-1)}e^{-2u^2+2uw}\RT),\ u\rw\IF,\ \delta\rw 0.
\EQNY

We have for $1\leq i< j\leq n $ and $u>0$
\BQNY
&&\pk{\sup_{\vk{x}\in E_i(\delta)}\LT(W(\vk{x})\Big|W(\vk{1})=w\RT)>u,
\sup_{\vk{x}\in E_j(\delta)}\LT(W(\vk{x})\Big|W(\vk{1})=w\RT)>u}\\
&&\leq \pk{\sup_{(\x,\y)\in E_i(\delta)\times E_j(\delta)}W(\x)+W(\y)-(F(\x)+F(\y))W(\vk{1})+(F(\x)+F(\y))w>2u}\\
&&\leq \pk{\sup_{(\x,\y)\in E_i(\delta)\times E_j(\delta)}W(\x)+W(\y)-(F(\x)+F(\y))W(\vk{1})>2u-2w}.
\EQNY
Since
\BQNY
\sigma^2_m&:=&\sup_{(\x,\y)\in E_i(\delta)\times E_j(\delta)}Var\LT(W(\x)+W(\y)-(F(\x)+F(\y))W(\vk{1})\RT)\\
&=&\sup_{(\x,\y)\in E_i(\delta)\times E_j(\delta)}\E{\LT(W(\x)+W(\y)-(F(\x)+F(\y))W(\vk{1})\RT)^2}\\
&=&\sup_{(\x,\y)\in E_i(\delta)\times E_j(\delta)}\LT(F(\x)+F(\y)\RT)\LT[1-\LT(F(\x)+F(\y)\RT)\RT]+2F(\x\wedge\y)\\
&<&\sup_{(\x,\y)\in E_i(\delta)\times E_j(\delta)}\LT(F(\x)+F(\y)\RT)\LT[2-\LT(F(\x)+F(\y)\RT)\RT]\\
&\leq& 1,
\EQNY
where we use the fact that for $(\x,\y)\in E_i(\delta)\times E_j(\delta)$
$$2F(\x\wedge\y)< F(\x)+F(\y).$$
By Borell-TIS inequality (ref.\cite{AdlerTaylor})
\BQNY
\pk{\sup_{(\x,\y)\in E_i(\delta)\times E_j(\delta)}W(\x)+W(\y)-(F(\x)+F(\y))W(\vk{1})>2u-2w}
&\leq& \exp\LT(-\frac{(2u-2w)^2}{2\sigma^2_m}\RT)\\
&=&o\LT( u^{2(n-1)}e^{-2u^2+2uw}\RT),\ u\rw\IF.
\EQNY
Thus we have
\BQNY
&&\sum_{1\leq i<j\leq n}\pk{\sup_{\vk{x}\in E_i(\delta)}\LT(W(\vk{x})\Big|W(\vk{1})=w\RT)>u,
\sup_{\vk{x}\in E_j(\delta)}\LT(W(\vk{x})\Big|W(\vk{1})=w\RT)>u}\\
&&\leq n^2\exp\LT(-\frac{(2u-2w)^2}{2\sigma^2_m}\RT)\\
&&=o\LT( u^{2(n-1)}e^{-2u^2+2uw}\RT),\ u\rw\IF.
\EQNY
Consequently, letting $u\rw\IF, \delta\rw 0$, we have
\BQNY
\pk{\sup_{\vk{x}\in [0,1]^n}\LT(W(\vk{x})\Big|W(\vk{1})=w\RT)>u}&\sim&
\sum_{i=1}^n\pk{\sup_{\vk{x}\in E_i(\delta)}\LT(W(\vk{x})\Big|W(\vk{1})=w\RT)>u}\\
&\sim& n(4d)^{(n-1)} \int_{L} \frac{\LT(\Pi_{i=1}^{n-1}x_i\RT)^{n-2}}{\LT(d+(1-d)\Pi_{i=1}^{n-1}x_i\RT)^{n-1}}d \wx u^{2(n-1)}e^{-2u^2+2uw}.
\EQNY
Specially, when $n=2$, we have $L=\LT[\frac{1}{\sqrt{(1-d)^2+1}+d},1\RT]$ and
$$\pk{\sup_{\vk{x}\in [0,1]^n}\LT(W(\vk{x})\Big|W(\vk{1})=w\RT)>u}\sim \frac{8d}{1-d}\ln \LT(\sqrt{1+\frac{1}{(1-d)^2}}-\frac{d}{1-d}\RT)u^{2}e^{-2u^2+2uw},\ u\rw\IF.$$
Further, if we take $\mathcal{E}=[0,1]^n$ and $\seEd:=\LT\{\x\in [0,1]^n: \frac{1}{2}-\delta \leq F(\vk{x})\leq \frac{1}{2}+\delta \RT\}$ with $\delta\in(0,\frac{1}{4})$ such that $d\LT(\frac{1}{2}-\delta\RT)+(1-d)\LT(\frac{1}{2}-\delta\RT)^n>\delta$,
then we have
$$\min_{1\leq i\leq n} x_i\geq \frac{1}{2}-\delta , \ \x \in \seEd,$$
and
\BQNY
 \inf_{\x,\y\in\seEd} F(\x\wedge\y)= \inf_{\x,\y\in\seEd} \LT(d\min_{1\leq i\leq n}x_i\wedge y_i+(1-d)\Pi_{i=1}^n(x_i\wedge y_i)\RT)
 \geq d\LT(\frac{1}{2}-\delta\RT)+(1-d)\LT(\frac{1}{2}-\delta\RT)^n>\delta,
\EQNY
which show that \eqref{Con2} holds. Thus by \nekorr{Corr3} and \netheo{Corr4}
\BQNY
&&\pk{\sup_{\vk{x}\in [0,1]^n}\abs{W(\vk{x})}>u\Big|W(\vk{1})=w}\nonumber\\
&&\sim \pk{\sup_{\vk{x}\in [0,1]^n}W(\vk{x})>u\Big|W(\vk{1})=w}+\pk{\sup_{\vk{x}\in [0,1]^n}W(\vk{x})>u\Big|W(\vk{1})=-w}\\
&&\sim c\pk{\sup_{\vk{x}\in [0,1]^n}\LT(W(\vk{x})\Big|W(\vk{1})=\abs{w}\RT)>u}
\EQNY
holds with $c=1$ for $w\neq 0$ and $c=2$ for $w=0$.

\QED

\proofprop{PROP1}
Clearly, for $F(\vk{x})$, \eqref{condition1}  holds for $\vk{x}, \vk{y}\in [0,1]^2$, hence the claim follows by  \netheo{Thm2}, \netheo{Corr4} and Remark \ref{Rem2}.

\QED

\COM{
\proofprop{PROP2}
For $\delta\in\LT(0,\frac{\sqrt{5}-2}{2}\RT)$ small enough, set
\BQNY
E_1(\delta)=\{\x\in [0,1]^2: x_1\leq x_2-\delta\},\
E_2(\delta)=\{\x\in [0,1]^2: x_1\geq x_2+\delta\},\
E_3(\delta)=\{\x\in [0,1]^2: x_2+\delta \leq x_1\leq x_2+\delta\}.
\EQNY
Then we have
\BQNY
\pk{\sup_{\vk{x}\in [0,1]^2}\LT(W(\vk{x})\Big|W(\vk{1})=w\RT)>u}&\leq&
\sum_{i=1}^3\pk{\sup_{\vk{x}\in E_i(\delta)}\LT(W(\vk{x})\Big|W(\vk{1})=w\RT)>u}\\
\pk{\sup_{\vk{x}\in [0,1]^2}\LT(W(\vk{x})\Big|W(\vk{1})=w\RT)>u}&\geq&
\sum_{i=1}^2\pk{\sup_{\vk{x}\in E_i(\delta)}\LT(W(\vk{x})\Big|W(\vk{1})=w\RT)>u}\\
&&-\pk{\sup_{\vk{x}\in E_1(\delta)}\LT(W(\vk{x})\Big|W(\vk{1})=w\RT)>u,
\sup_{\vk{x}\in E_2(\delta)}\LT(W(\vk{x})\Big|W(\vk{1})=w\RT)>u}.
\EQNY
For $\x\in E_1(\delta)$, we have $F(\x)=\frac{1}{2}x_1\LT(1+x_2\RT)$,
\BQNY
\seE_1=\{\x\in E_1(\delta):F(\x)=\frac{1}{2}\}=\LT\{(x_1,h(x_1)):x_1\in
\LT[\frac{1}{2},\frac{\sqrt{5}-1}{2}-\delta\RT]\RT\},\ h(x_1)=\frac{1-x_1}{x_1},
\EQNY
and
\BQNY
\lim_{\vn\rw 0}\sup_{\vk{z}\in \seE_1 }\underset{\vk{x}\neq \vk{y}}{\sup_{\abs{\vk{x}-\vk{z}},\abs{\vk{y}-\vk{z}}\leq \vn}}\frac{\abs{F(\vk{x})-F(\vk{y})-\frac{1}{2}(1+z_2)(x_1-y_1)-\frac{1}{2}z_1(x_2-y_2)}}
{\abs{x_1-y_1}+\abs{x_2-y_2}}=0,
\EQNY
Thus by \netheo{Thm1}, we have
\BQNY
\pk{\sup_{\vk{x}\in E_1(\delta)}\LT(W(\vk{x})\Big|W(\vk{1})=w\RT)>u}
\sim 8\int_{\frac{1}{2}}^{\frac{\sqrt{5}-1}{2}-\delta} \frac{1}{2x}d x u^2e^{-2u^2+2uw}
=4\ln (\sqrt{5}-1-2\delta) u^2e^{-2u^2+2uw}, u\rw\IF.
\EQNY
By the symmetry of $F(\x)$ on $E_1(\delta)$ and $E_1(\delta)$,
\BQNY
\pk{\sup_{\vk{x}\in E_2(\delta)}\LT(W(\vk{x})\Big|W(\vk{1})=w\RT)>u}
\sim 4\ln (\sqrt{5}-1-2\delta) u^2e^{-2u^2+2uw}, u\rw\IF.
\EQNY
Assume that  $W_1(\x)$ is a Brownian sheet based on $F_1(\x)=\frac{1}{2}x_1(1+x_2)$ and
$W_2(\x)$ is a Brownian sheet based on $F_2(\x)=\frac{1}{2}x_2(1+x_1)$. Then
\BQNY
\pk{\sup_{\vk{x}\in E_3(\delta)}\LT(W(\vk{x})\Big|W(\vk{1})=w\RT)>u}
\leq \pk{\sup_{\vk{x}\in E_3(\delta)}\LT(W_1(\vk{x})\Big|W_1(\vk{1})=w\RT)>u}
+\pk{\sup_{\vk{x}\in E_3(\delta)}\LT(W(\vk{x})\Big|W(\vk{1})=w\RT)>u}
\EQNY
For $F_1(\x)=\frac{1}{2}x_1(1+x_2)$ and $\x\in E_3(\delta)$ we have
\BQNY
\seE_2=\{\x\in E_1(\delta):F_1(\x)=\frac{1}{2}\}=\LT\{(x_1,h(x_1)):x_1\in
\LT[\frac{\sqrt{5}-1}{2}-\delta,\frac{\sqrt{5}-1}{2}+\delta \RT]\RT\},\ h(x_1)=\frac{1-x_1}{x_1},
\EQNY
and
\BQNY
\lim_{\vn\rw 0}\sup_{\vk{z}\in \seE_2 }\underset{\vk{x}\neq \vk{y}}{\sup_{\abs{\vk{x}-\vk{z}},\abs{\vk{y}-\vk{z}}\leq \vn}}\frac{\abs{F_1(\vk{x})-F_1(\vk{y})-\frac{1}{2}(1+z_2)(x_1-y_1)-\frac{1}{2}z_1(x_2-y_2)}}
{\abs{x_1-y_1}+\abs{x_2-y_2}}=0.
\EQNY
Further by \netheo{Thm1}, we have
\BQNY
\pk{\sup_{\vk{x}\in E_3(\delta)}\LT(W_1(\vk{x})\Big|W_1(\vk{1})=w\RT)>u}
\sim 4\ln \frac{\sqrt{5}-1+2\delta}{\sqrt{5}-1-2\delta} u^2e^{-2u^2+2uw}
, u\rw\IF.
\EQNY
Similarly,
\BQNY
\pk{\sup_{\vk{x}\in E_3(\delta)}\LT(W_2(\vk{x})\Big|W_2(\vk{1})=w\RT)>u}
\sim 4\ln \frac{\sqrt{5}-1+2\delta}{\sqrt{5}-1-2\delta} u^2e^{-2u^2+2uw}
, u\rw\IF.
\EQNY
We have for $u>0$
\BQNY
&&\pk{\sup_{\vk{x}\in E_1(\delta)}\LT(W(\vk{x})\Big|W(\vk{1})=w\RT)>u,
\sup_{\vk{x}\in E_2(\delta)}\LT(W(\vk{x})\Big|W(\vk{1})=w\RT)>u}\\
&&\leq \pk{\sup_{(\x,\y)\in E_1(\delta)\times E_2(\delta)}W(\x)+W(\y)-(F(\x)+F(\y))W(\vk{1})+(F(\x)+F(\y))w>2u}\\
&&\leq \pk{\sup_{(\x,\y)\in E_1(\delta)\times E_2(\delta)}W(\x)+W(\y)-(F(\x)+F(\y))W(\vk{1})>2u-2w}.
\EQNY
Since
\BQNY
\sigma^2_m&:=&\sup_{(\x,\y)\in E_1(\delta)\times E_2(\delta)}Var\LT(W(\x)+W(\y)-(F(\x)+F(\y))W(\vk{1})\RT)\\
&=&\sup_{(\x,\y)\in E_1(\delta)\times E_2(\delta)}\E{\LT(W(\x)+W(\y)-(F(\x)+F(\y))W(\vk{1})\RT)^2}\\
&=&\sup_{(\x,\y)\in E_1(\delta)\times E_2(\delta)}\LT(F(\x)+F(\y)\RT)\LT[1-\LT(F(\x)+F(\y)\RT)\RT]+2F(\x\wedge\y)\\
&<&\sup_{(\x,\y)\in E_1(\delta)\times E_2(\delta)}\LT(F(\x)+F(\y)\RT)\LT[2-\LT(F(\x)+F(\y)\RT)\RT]\\
&\leq& 1,
\EQNY
where we use the fact that for $(\x,\y)\in E_1(\delta)\times E_2(\delta)$
$$2F(\x\wedge\y)< F(\x)+F(\y).$$
Then by Borell inequality
\BQNY
\pk{\sup_{(\x,\y)\in E_1(\delta)\times E_2(\delta)}W(\x)+W(\y)-(F(\x)+F(\y))W(\vk{1})>2u-2w}
&\leq& \exp\LT(-\frac{(2u-2w)^2}{2\sigma^2_m}\RT)\\
&=&o\LT( u^2e^{-2u^2+2uw}\RT),\ u\rw\IF.
\EQNY

Consequently, letting $u\rw\IF, \delta\rw 0$, we have
\BQNY
\pk{\sup_{\vk{x}\in [0,1]^2}\LT(W(\vk{x})\Big|W(\vk{1})=w\RT)>u}\sim
\sum_{i=1}^2\pk{\sup_{\vk{x}\in E_i(\delta)}\LT(W(\vk{x})\Big|W(\vk{1})=w\RT)>u}
\sim 8\ln (\sqrt{5}-1) u^2e^{-2u^2+2uw}.
\EQNY

\QED
}

\def\wk{\widetilde{\vk{k}}}
\def\k{\vk{k}}
\def\wl{\widetilde{\vk{l}}}
\def\l{\vk{l}}
\section{Appendix}
Before stating the proofs of next lemmas, we introduce some notation.
Define the random fields $Y(\x), \x \in \R^n$ by
\BQN\label{Yc}
Y(\vk{x})=\sum_{i=1}^{n}B^{(i)}(x_i)
\EQN
where $B^{(i)}$'s are independent  standard Brownian motions. We define
\BQNY
\mathcal{H}_{\vk{c}} [\vk{0},\vk{\lambda}]=\E{\sup_{\x\in[0,\vk{c}*\vk{\lambda}]}e^{Y(\x)-\Var(Y (\vk{x}))}}
\EQNY
and notice that
\BQN\label{HH2}
\lim_{\min_{1 \le i \le n} {\lambda_i}\rw\vk{\IF}}
\frac{1}{\Pi_{i=1}^n\lambda_i}\mathcal{H}_{\vk{c}} [\vk{0},\vk{\lambda}]
= \E{ \frac{\sup_{\vk x\in \R^n }e^{ \sqrt{2}Y (\vk{x})-\Var(Y (\vk{x}))}}{\int_{\R^n }e^{ \sqrt{2}Y (\vk{x})-\Var(Y (\vk{x})) } d \vk x  }} = \Pi_{i=1}^nc_i .
\EQN
See the recent contributions \cite{EHKD,SBK} for various results on Pickands constants.
\BL\label{lem1}
Let $h(\wx),\wx=(x_1\ldot x_{n-1})\in\R^{n-1}$ is a continuous differentiable function. Assume that $X(\vk{x}),\ \vk{x}\in E,$ with $E=\{\vk{x}: x_i\in[S_i,T_i], i=1\ldot   n-1,\abs{x_n-h(\wx)}\leq T\}, T>0, 0\leq S_i<T_i, i=1\ldot {n-1}$ is a Gaussian field which has continuous pathes, variance function $\sigma^2(\vk{x})$ and correlation function $r(\vk{x}, \vk{y})$ and $g(\vk{x}), \vk{x}\in E$ is a continuous function. Further, $\sigma^2(\vk{x})$ attains it maximum at $\widetilde{E}=\{\vk{x}: x_i\in[S_i,T_i], i=1\ldot   n-1, x_n=h(\wx)\}$ which satisfies for $\x\in E$
\BQN\label{var2}
\lim_{\delta\rw 0}\underset{\x\in E\setminus \widetilde{E}}{\sup_{\abs{\x-h(\wx)}\leq\delta}}
\abs{\frac{1-\sigma(\vk{x})}{b\LT(x_n-h(\wx)\RT)^2}-1}=0,
\EQN
and
\BQN\label{g2}
\lim_{\delta\rw 0}\underset{\x\in E\setminus \widetilde{E}}{\sup_{\abs{\x-h(\wx)}\leq\delta}}
\abs{\frac{g(\x)}{c\LT(x_n-h(\wx)\RT)}-1}=0,
\EQN
for constants $b>0$ and $c\in\R$.\\
For any $\z\in\widetilde{E}$, $\vk{x},\ \vk{y}\in E$
\BQN\label{r2}
1-r(\vk{x},\vk{y})\sim \sum_{i=1}^{n}c_i\abs{x_i-y_i}, \ \abs{\x-\z}, \ \abs{\y-\z} \rw0,
\EQN
 where $ c_i>0$ are constants and there exist positive constants $\mathcal{C}_1$ and $\mathcal{C}_2$ such that
\BQN\label{r3}
\mathcal{C}_1\sum_{i=1}^{n}\abs{x_i-y_i}<1-r(\vk{x},\vk{y})< \mathcal{C}_2\sum_{i=1}^{n}\abs{x_i-y_i} , \ \vk{x},\ \vk{y}\in E.
\EQN
If further there exists a positive constant $\mathcal{C}_3$ such that
\BQN\label{boundh}
\underset{\x\neq\y}{\sup_{\x,\y\in E}}\frac{\abs{(x_n-h(\wx))-(y_n-h(\wy))}}{\sum_{i=1}^{n} \abs{x_i-y_i}}
\leq \mathcal{C}_3 .
\EQN
\COM{\BQN\label{boundv}
\underset{\x\neq\y}{\sup_{\x,\y\in E}}\frac{\abs{\frac{1}{\sigma(\x)}-\frac{1}{\sigma(\y)}}}{(\abs{x_n-h(\wx)}+\abs{y_n-h(\wy)})\sum_{i=1}^{n}\abs{x_i-y_i}}\leq \mathcal{C}_3,
\EQN
and
\BQN\label{boundg}
\underset{\x\neq\y}{\sup_{\x,\y\in E}}\frac{\abs{g(\x)-g(\y)}}{\sum_{i=1}^n\abs{x_i-y_i}}\leq \mathcal{C}_3,
\EQN}
then as $u\rw\IF$
\BQNY
\pk{\sup_{\vk{x}\in E}(X(\vk{x})+g(\vk{x}))>u}&\sim&
\LT(\Pi_{i=1}^{n-1}(T_i-S_i)\RT)\LT(\Pi_{i=1}^nc_i\RT) \int_{-\IF}^{\IF}e^{-bx^{2}+cx}dx u^{2n-1}\Psi(u)\\
&=&\LT(\Pi_{i=1}^{n-1}(T_i-S_i)\RT)\LT(\Pi_{i=1}^nc_i\RT)
\sqrt{\frac{\pi}{b}}e^{\frac{c^2}{4b}}u^{2n-1}\Psi(u).
\EQNY
\EL
\prooflem{lem1}
In the following proof, without loss of generality, we assume that $c>0$.\\
First for $\vn_1,\vn\in(0,1)$ and $\lambda>0$ we introduce the following notation:
\begin{align*}
&E_0=\Pi_{i=1}^{n-1}[S_i,T_i],\ \wx=(x_1\ldot  x_{n-1}),\
E(\vn_1)=\LT\{\vk{x}\in E: \abs{x_n-h(\wx)}\leq \vn_1,\ \wx\in E_0\RT\},\\
& E(u)=\LT\{\vk{x}:  \abs{x_n-h(\wx)}\leq \frac{\ln u}{u},\  \wx\in E_0\RT\},\
J_{k}(u)=\LT[\frac{k\lambda}{u^{2}},\frac{(k+1)\lambda}{u^{2}}\RT], \
 M_{n}(u)=\LT\lfloor\frac{u\ln u}{\lambda}\RT\rfloor,\\
&\wk=(k_1\ldot k_{n-1}),\ D_{\vk{k}}(u)=\Pi_{i=1}^nJ_{k_{i}}(u),\ \k\in\mathbb{N}^n,\
 D_{\wk}(u)=\Pi_{i=1}^{n-1}J_{k_{i}}(u),\ \wk\in\mathbb{N}^{n-1},\\
& \mathcal{M}_1(u)=\LT\{\wk:D_{\wk}(u)\subset \prod_{i=1}^{n-1}\LT[S_i,T_i\RT] \RT\},\
\mathcal{M}_2(u)=\LT\{\wk:D_{\wk}(u)\cap \prod_{i=1}^{n-1}\LT[S_i,T_i\RT]\neq \emptyset \RT\},\\
&\mathcal{L}_1(u)=\{\vk{k}: D_{\vk{k}}(u)\subset E(u)\},\
\mathcal{L}_2(u)=\{\vk{k}: D_{\vk{k}}(u)\cap E(u)\neq \emptyset\}.\\
&\mathcal{K}_1(u)=\{(\vk{k},\vk{l}):\vk{k},\vk{l}\in\mathcal{L}_1(u),\vk{k}\neq \vk{l}, D_{\vk{k}}(u)\cap D_{\vk{l}}(u)\neq \emptyset\},\\
&\mathcal{K}_2(u)=\{(\vk{k},\vk{l}):\vk{k},\vk{l}\in\mathcal{L}_1(u),D_{\vk{k}}(u)\cap D_{\vk{l}}(u)=\emptyset, u^{-2}\abs{k_1-l_1}\lambda\leq \vn\},\\
&\mathcal{K}_3(u)=\{(\vk{k},\vk{l}):\vk{k},\vk{l}\in\mathcal{L}_1(u),D_{\vk{k}}(u)\cap D_{\vk{l}}(u)=\emptyset, u^{-2}\abs{k_1-l_1}\lambda\geq \vn\}.
\end{align*}
Here we need to notice that $\vk{i}, \vk{k}$ and $\vk{l}$ are $n$-dimensional.\\
It follows that for $u$ \large enough
\BQN\label{sli1}
\Pi_0(u)\leq \pk{\sup_{\vk{x}\in E}(X(\vk{x})+g(\vk{x}))>u}\leq \Pi_0(u)+\Pi_1(u)+\Pi_2(u)
\EQN
where
\BQNY
&&\Pi_0(u):=\pk{\sup_{\vk{x}\in E(u)}(X(\vk{x})+g(\vk{x}))>u},\quad
\Pi_1(u):=\pk{\sup_{\vk{x}\in E\setminus E(\vn_1)}(X(\vk{x})+g(\vk{x}))>u},\\
&&\Pi_2(u):=\pk{\sup_{\vk{x}\in E(\vn_1)\setminus E(u)}(X(\vk{x})+g(\vk{x}))>u}.
\EQNY
By \eqref{var2} and \eqref{g2}, we have for any small $\vn\in(0,1)$ there exists $\vn_1$ small enough such that
\BQN\label{var3}
1-(1+\vn)b(x_n-h(\wx))^2\leq \sigma(x)\leq 1-(1-\vn)b(x_n-h(\wx))^2
\EQN
and
\BQN\label{g3}
c(x_n-h(\wx))-\vn\abs{x_n-h(\wx)}\leq g(\wx)\leq c(x_n-h(\wx))+\vn\abs{x_n-h(\wx)}
\EQN
hold for $\x\in E(\vn_1)$.\\
By continuity of $\sigma(\x)$ and $ E(\vn_1)\supset \widetilde{E}$, we have
$$\sup_{\vk{x}\in E\setminus E(\vn_1)}\sigma(\vk{x})<1-\delta_1$$
with $\delta_1\in (0,1)$, which combined with Borell-TIS inequality as in \cite{AdlerTaylor} leads
\BQN\label{boundofpi1}
\Pi_1(u)\leq e^{-\frac{\LT(u-\mathbb{Q}_1-\mathbb{Q}_2\RT)^2}{2(1-\delta_1)^2}}
=o\LT(\Psi(u)\RT),\ u\rw\IF,
\EQN
where $\mathbb{Q}_1=\sup_{\vk{x}\in E\setminus E(\vn_1)}g(\vk{x})<\IF$ and $\mathbb{Q}_2=\E{\sup_{\vk{x}\in E\setminus E(\vn_1)} X(\vk{x})}<\IF$.\\
In light of \eqref{var3} and \eqref{g3}, we have for $u$ large enough
\BQNY
\inf_{\vk{x}\in E(\vn_1)\setminus E(u)}\frac{1}{\sigma(\vk{x})}&\geq& 1+\mathbb{Q}_3\LT(\frac{\ln u}{u}\RT)^2,
\EQNY
and
\BQNY
\sup_{\vk{x}\in E(\vn_1)\setminus E(u)}g(\vk{x})&\leq&\mathbb{Q}_4\LT(\frac{\ln u}{u}\RT).
\EQNY
By \eqref{r3}, we have
\BQNY
\E{\LT(\overline{X}(\vk{x})-\overline{X}(\vk{y})\RT)^2}
=2(1-r(\vk{x},\vk{y}))\leq 2\mathcal{C}_2\sum_{i=1}^{n}\abs{x_i-y_i},
\EQNY
which combined with  \cite{Pit96} [Theorm8.1] derive for $u$ large enough
\BQN\label{err0}
 \Pi_2(u) &\leq&
\pk{\sup_{\vk{x}\in E(\vn_1)\setminus E(u)}X(\vk{x})>u-\mathbb{Q}_4\LT(\frac{\ln u}{u}\RT)}\nonumber\\
&\leq&
\pk{\sup_{\vk{x}\in E(\vn_1)\setminus E(u)}\overline{X}(\vk{x})>\LT(u-\mathbb{Q}_4\LT(\frac{\ln u}{u}\RT)\RT)\LT(1+\mathbb{Q}_3\LT(\frac{\ln u}{u}\RT)^2\RT)}\nonumber\\
&\leq&\mathbb{Q}_5u^{2n}\Psi\LT(\LT(u-\mathbb{Q}_4\LT(\frac{\ln u}{u}\RT)\RT)\LT(1+\mathbb{Q}_3\LT(\frac{\ln u}{u}\RT)^2\RT)\RT)\nonumber\\
&=&o\LT(\Psi(u)\RT),\ u\rw\IF,
\EQN
By \eqref{sli1}, \eqref{boundofpi1}, \eqref{err0} and the fact $ \Pi_0(u)\geq \pk{X(\z)>u}=\Psi(u)$ with $\z\in\widetilde{E}$ leads to
\BQN\label{asym1}
\pk{\sup_{\vk{x}\in E}(X(\vk{x})+g(\vk{x}))>u}\sim \Pi_0(u), u\rw\IF.
\EQN
Next we focus on $\Pi_0(u)$.
By \eqref{boundh}, we have
\BQNY
\sup_{\vk{k}\in\mathcal{L}_2(u)}\sup_{\x,\y\in D_{\vk{k}}(u)}
\abs{(x_n-h(\wx))-(y_n-h(\wy))}\leq \mathbb{Q}_6 \frac{\lambda}{u^2}
\EQNY
and
\BQNY
&&\sup_{\vk{k}\in\mathcal{L}_2(u)}\sup_{\x,\y\in D_{\vk{k}}(u)}\abs{(x_n-h(\wx))^2-(y_n-h(\wy))^2}\\
&&\leq
\sup_{\vk{k}\in\mathcal{L}_2(u)}\sup_{\x,\y\in D_{\vk{k}}(u)}(\abs{x_n-h(\wx)}+\abs{y_n-h(\wy)})\abs{(x_n-h(\wx))-(y_n-h(\wy))}\\
&&\leq\mathbb{Q}_7 \frac{\lambda\ln u }{u^3}.
\EQNY
In the view of  \eqref{var3} and \eqref{g3} we notice that for any $\vk{k}\in\mathcal{L}_2(u)$
\BQNY
u_{\vk{k}}^{+}&:=&\sup_{\x\in D_{\vk{k}}(u)}(u-g(\x))\frac{1}{\sigma(\x)}\\
&\leq& \sup_{\x\in D_{\vk{k}}(u)}(u-c(x_n-h(\wx))+\vn \abs{x_n-h(\wx)})
(1+(1+\vn)b(x_n-h(\wx))^2)\\
&\leq&\LT(u-\inf_{\x\in D_{\vk{k}}(u)}(c(x_n-h(\wx)))+\vn\sup_{\x\in D_{\vk{k}}(u)}\abs{x_n-h(\wx)}\RT)\LT(1+(1+\vn)b\sup_{\x\in D_{\vk{k}}(u)}(x_n-h(\wx))^2\RT)\\
&\leq& \LT(u-\inf_{y\in J_{l}(u)}(cy)+\vn\sup_{y\in  J_{l}(u) }\abs{y}+\mathbb{Q}_8 \frac{\lambda}{u^2}\RT)\LT(1+(1+\vn)b\sup_{y\in J_{l}(u)}(y)^2+\mathbb{Q}_9 \frac{\lambda\ln u }{u^3}\RT)\\
&=:& u^{+}_l,
\EQNY
where $l$ satisfy that
\BQN\label{conditionl}
J_{l}(u)\cap \LT[\inf_{\x\in D_{\vk{k}}(u)}(x_n-h(\wx)),\sup_{\x\in D_{\vk{k}}(u)}(x_n-h(\wx))\RT]\neq \emptyset.
\EQN
Similarly, we define
\BQNY
u_{\vk{k}}^{-}:=\inf_{\x\in D_{\vk{k}}(u)}(u-g(\x))\frac{1}{\sigma(\x)}
\EQNY
and
\BQNY
u^{-}_l= \LT(u-\sup_{y\in J_{l}(u)}(cy)-\vn\sup_{y\in  J_{l}(u) }\abs{y}-\mathbb{Q}_8 \frac{\lambda}{u^2}\RT)\LT(1+(1-\vn)b\inf_{y\in J_{l}(u)}(y)^2-\mathbb{Q}_9 \frac{\lambda\ln u }{u^3}\RT)
\EQNY
where $u_{\vk{k}}^{-}\geq u^{-}_l$ with $l$ satisfying \eqref{conditionl}.

Considering $\k\in\mathcal{L}_2(u)$, if we fix $\wk$ first, for all $k_n$ such that
$\k\in\mathcal{L}_2(u)$ we can chose $l$ satisfying \eqref{conditionl}  from $-M_n(u)-1$ to $M_n(u)+1$.\\
Bonferroni inequality leads to
\BQN
&&\Pi_0(u)\leq\sum_{\vk{k}\in \mathcal{L}_2(u)}
\pk{\sup_{\vk{x}\in D_{\vk{k}}(u)}(X(\vk{x})+g(\vk{x}))>u},\label{up1}\\
&&\Pi_0(u)\geq\sum_{\vk{k}\in \mathcal{L}_1(u)}
\pk{\sup_{\vk{x}\in D_{\vk{k}}(u)}(X(\vk{x})+g(\vk{x}))>u}-\sum_{i=1}^{3}\mathcal{A}_i(u),\label{low1}
\EQN
where
\BQNY
\mathcal{A}_i(u)&=&\sum_{(\vk{k}, \vk{l})\in\mathcal{K}_i(u)}\pk{\sup_{\vk{x}\in D_{\vk{k}}(u)}(X(\vk{x})+g(\vk{x}))>u,
\sup_{\vk{x}\in D_{\vk{l}}(u)}(X(\vk{x})+g(\vk{x}))>u}\\
&\leq&\sum_{(\vk{k},\vk{l})\in\mathcal{K}_i(u)}\pk{\sup_{\vk{x}\in D_{\vk{k}}(u)}\overline{X}(\vk{x})>u_{\vk{k}}^{-},
\sup_{\vk{x}\in D_{\vk{l}}(u)}\overline{X}(\vk{x})>u_{\vk{l}}^{-}}, i=1,2,3.
\EQNY
We set
\BQNY
X_{u,\vk{k}}(\vk{x})=\overline{X}(k_1u^{-2}\lambda+x_1\ldot k_nu^{-2}\lambda+x_n),\ \vk{x}\in D_{\vk{0}}(u), \vk{k}\in\mathcal{L}_2(u).
\EQNY
Then by \eqref{r2} and \nelem{lem0} that
\BQN\label{EQA11}
\lim_{u\rw\IF}\sup_{\vk{k}\in\mathcal{L}_2(u)}\abs{
\frac{\pk{\sup_{\vk{x}\in D_{\vk{0}}(u)}X_{u,\vk{k}}(\vk{x})>u_{\vk{k}}^{-}}}{\Psi(u_{\vk{k}}^{-})}
-\mathcal{H} _{\vk{c}}[\vk{0},\vk{\lambda}]
}=0,
\EQN
where $\vk{\lambda}=(\lambda\ldot  \lambda).$
Further, by \eqref{HH2}
\begin{align}\label{upper}
&\sum_{\vk{k}\in \mathcal{L}_2(u)}
\pk{\sup_{\vk{x}\in D_{\vk{k}}(u)}(X(\vk{x})+g(\vk{x}))>u}\nonumber\\
&\leq
\sum_{\vk{k}\in \mathcal{L}_2(u)}\pk{\sup_{\vk{x}\in D_{\vk{0}}(u)}X_{u,\vk{k}}(\vk{x})>u_{\vk{k}}^{-}}\nonumber\\
&\sim \mathcal{H}_{\vk{c}}[\vk{0},\vk{\lambda}]
\sum_{\vk{k}\in \mathcal{L}_2(u)}\Psi(u_{\vk{k}}^{-})\nonumber\\
&\leq \mathcal{H} _{\vk{c}}[\vk{0},\vk{\lambda}]
\sum_{\wk\in\mathcal{M}_2(u)}\LT(\sum_{l=-M_{n}(u)-1}^{M_{n}(u)+1}
\Psi(u_{l}^{-})\RT)\nonumber\\
&\sim \mathcal{H} _{\vk{c}}[\vk{0},\vk{\lambda}]\Psi(u)
\sum_{\wk\in\mathcal{M}_2(u)}\LT(\sum_{l=-M_{n}(u)-1}^{M_{n}(u)+1}
e^{\sup_{y\in [l,l+1]}\LT(cy+\vn\abs{y}\RT)\frac{\lambda}{u}-(1-\vn)b\inf_{y\in [l,l+1]}\frac{(y\lambda)^2}{u^2}}\RT)\nonumber\\
&\sim \mathcal{H} _{\vk{c}}[\vk{0},\vk{\lambda}]\Psi(u)
\sum_{\wk\in\mathcal{M}_2(u)}\frac{u}{\lambda}\int_{-\IF}^{\IF}
e^{-(1-\vn)bx^2+cx+\vn\abs{x}}dx\nonumber\\
&\sim \frac{\mathcal{H} _{\vk{c}}[\vk{0},\vk{\lambda}]}{\lambda^n}\Psi(u)
\LT(\Pi_{i=1}^{n-1}(T_i-S_i)\RT)u^{2n-1}\int_{-\IF}^{\IF}e^{-(1-\vn)bx^2+cx+\vn\abs{x}}dx \nonumber\\
&\sim \LT(\Pi_{i=1}^{n-1}(T_i-S_i)\RT)\LT(\Pi_{i=1}^nc_i\RT) \int_{-\IF}^{\IF}e^{-bx^{2}+cx}dx u^{2n-1}\Psi(u),
\end{align}
as $u\rw\IF, \lambda\rw\IF, \vn\rw0$. Similarly,
\BQN
&&\sum_{\vk{k}\in \mathcal{L}_1(u)}
\pk{\sup_{\vk{x}\in D_{k,l}(u)}(X(\vk{x})+g(\vk{x}))>u}\geq
\sum_{\vk{k}\in \mathcal{L}_1(u)}\pk{\sup_{\vk{x}\in D_{\vk{0}}(u)}X_{u,\vk{k}}(\vk{x})>u_{\vk{k}}^{+}}\nonumber\\
&&\quad\quad\geq  \LT(\Pi_{i=1}^{n-1}(T_i-S_i)\RT)\LT(\Pi_{i=1}^nc_i\RT) \int_{-\IF}^{\IF}e^{-bx^{2}+cx}dx u^{2n-1}\Psi(u), \ u\rw\IF, \lambda\rw\IF, \vn\rw0.
\EQN
Next we will show that $\mathcal{A}_i(u), i=1,2,3$ are all negligible compared with
$$\sum_{\vk{k}\in \mathcal{L}_1(u)}\pk{\sup_{\vk{x}\in D_{\vk{k}}(u)}(X(\vk{x})+g(\vk{x}))>u}.$$
For any $(\vk{k},\vk{l})\in \mathcal{K}_1(u)$, without loss of generality, we assume that
$k_1+1=l_1$. Let
$$D_{\vk{k}}^1(u)=\LT[k\frac{\lambda}{u^{2}},((k+1)\lambda-\sqrt{\lambda})\frac{1}{u^{2}}\RT]
\times \Pi_{j=2}^{n}J_{k_j}(u), \ D_{\vk{k}}^2(u)=\LT[((k+1)\lambda-\sqrt{\lambda})\frac{1}{u^{2}},(k+1)\frac{\lambda}{u^{2}},\RT]
\times \Pi_{j=2}^{n}J_{k_j}(u).$$
Then
\BQNY
&&\pk{\sup_{\vk{x}\in D_{\vk{k}}(u)}\overline{X}(\vk{x})>u_{\vk{k}}^{-\vn},
\sup_{\vk{x}\in D_{\vk{l}}(u)}\overline{X}(\vk{x})>u_{\vk{l}}^{-\vn}}\\
&&\leq\pk{\sup_{\vk{x}\in D^1_{\vk{k}}(u)}\overline{X}(\vk{x})>u_{\vk{k}}^{-\vn},
\sup_{\vk{x}\in D_{\vk{l}}(u)}\overline{X}(\vk{x})>u_{\vk{l}}^{-\vn}}+\pk{\sup_{\vk{x}\in D^2_{\vk{k}}(u)}\overline{X}(\vk{x})>u_{\vk{k}}^{-\vn}}.
\EQNY
Analogously as in \eqref{EQA11}, we have
\BQNY
\lim_{u\rw\IF}\sup_{\vk{k}\in\mathcal{L}_1(u)}\abs{
\frac{\pk{\sup_{\vk{x}\in D^2_{\vk{k}}(u)}\overline{X}(\vk{x})>u_{\vk{k}}^{-\vn}}}{\Psi(u_{\vk{k}}^{-\vn})}
-\mathcal{H} _{\vk{c}}[\vk{0},\vk{\lambda}_1]
}=0,
\EQNY
where $\vk{\lambda}_1=(\sqrt{\lambda},\lambda\ldot  \lambda).$
Moreover, in the light of \eqref{r2} and \cite{EGRFRV2017} [Lemma 5.4] we have for $u$ large enough
\BQNY
\pk{\sup_{\vk{x}\in D^1_{\vk{k}}(u)}\overline{X}(\vk{x})>u_{\vk{k}}^{-\vn},
\sup_{\vk{x}\in D_{\vk{l}}(u)}\overline{X}(\vk{x})>u_{\vk{l}}^{-\vn}}\leq \mathbb{Q}_7\lambda^4e^{-\mathbb{Q}_8\lambda^{1/4}}
\Psi(\min(u_{\vk{k}}^{-\vn},u_{\vk{l}}^{-\vn})),
\EQNY
and for $(\vk{k},\vk{l})\in \mathcal{K}_2(u)$,
\BQNY
\pk{\sup_{\vk{x}\in D_{\vk{k}}(u)}\overline{X}(\vk{x})>u_{\vk{k}}^{-\vn},
\sup_{\vk{x}\in D_{\vk{l}}(u)}\overline{X}(\vk{x})>u_{\vk{l}}^{-\vn}}
\leq\mathbb{Q}_9\lambda^4e^{-\mathbb{Q}_{10}
\abs{\sum_{i=1}^n(k_i-l_i)^2}^{1/4}\lambda^{1/2}}
\Psi(\min(u_{\vk{k}}^{-\vn},u_{\vk{l}}^{-\vn})),
\EQNY
where $\mathbb{Q}_i, i=7,8,9,10$ are positive constants independent of $u$ and $\lambda$.\\
Since $D_{\vk{k}}(u)$ has at most $3^n-1$ neighbors, then
\BQN
\mathcal{A}_1(u)&\leq& \sum_{(\vk{k},\vk{l})\in\mathcal{K}_1(u)}\pk{\sup_{\vk{x}\in D_{\vk{k}}(u)}\overline{X}(\vk{x})>u_{\vk{k}}^{-\vn},
\sup_{\vk{x}\in D_{\vk{l}}(u)}\overline{X}(\vk{x})>u_{\vk{l}}^{-\vn}}\nonumber\\
&\leq& 2\sum_{(\vk{k},\vk{l})\in\mathcal{K}_1(u)}
\LT(\mathcal{H} _{\vk{c}}[\vk{0},\vk{\lambda}_1]
+\mathbb{Q}_7\lambda^4e^{-\mathbb{Q}_8\lambda^{1/4}}\RT)
\Psi(\min(u_{\vk{k}}^{-\vn},u_{\vk{l}}^{-\vn}))\nonumber\\
&\leq& 2\times (3^n-1)\sum_{\vk{k}\in\mathcal{L}_1(u)}
\LT(\mathcal{H} _{\vk{c}}[\vk{0},\vk{\lambda}_1]
+\mathbb{Q}_7\lambda^4e^{-\mathbb{Q}_8\lambda^{1/4}}\RT)
\Psi(u_{\vk{k}}^{-\vn})\nonumber\\
&=&o\LT(u^{2n-1}\Psi(u)\RT),\ u\rw\IF,\ \lambda\rw\IF.
\EQN
and
\BQN
\mathcal{A}_2(u)&\leq& \sum_{(\vk{k},\vk{l})\in\mathcal{K}_2(u)}\mathbb{Q}_9\lambda^4e^{-\mathbb{Q}_{10}
\abs{\vk{k}-\vk{l}}^{1/2}\lambda^{1/2}}
\Psi(\min(u_{\vk{k}}^{-\vn},u_{\vk{l}}^{-\vn}))\nonumber\\
&\leq&\sum_{\vk{k}\in\mathcal{L}_1(u)}\mathbb{Q}_9\lambda^4\Psi\LT(u_{\vk{k}}^{-\vn}\RT)
\sum_{\abs{\vk{k}}\geq 1 }e^{-\mathbb{Q}_{10}
\abs{\vk{k}}^{1/2}\lambda^{1/2}}\nonumber\\
&\leq&\sum_{\vk{k}\in\mathcal{L}_1(u)}\mathbb{Q}_9\lambda^4\Psi\LT(u_{\vk{k}}^{-\vn}\RT)
e^{-\mathbb{Q}_{10}\lambda^{1/2}}\nonumber\\
&=&o\LT(u^{2n-1}\Psi(u)\RT), \ u\rw\IF, \lambda\rw\IF.
\EQN
For $(\vk{k},\vk{l})\in \mathcal{K}_3(u)$, $\abs{x_{n}-y_{n}}\geq \vn/2$ holds with
$\vk{x}\in D_{\vk{k}}(u), \ \vk{y}\in D_{\vk{l}}(u)$. Then by \eqref{r3}, for $u$ large enough
\BQNY
\Var\LT(\overline{X}(\vk{x})+\overline{X}(\vk{y})\RT)=2(1+r(\vk{x},\vk{y}))
\leq2+2\sup_{\abs{x_1-y_1}\geq \vn/2}r(\vk{x},\vk{y})\leq 4-\delta
\EQNY
holds with $\delta\in(0,1)$ for $(\vk{k},\vk{l})\in \mathcal{K}_3(u),
\vk{x}\in D_{\vk{k}}(u), \vk{y}\in D_{\vk{l}}(u).$ Further, Borell-TIS inequality leads to
\BQN\label{err3}
\mathcal{A}_3(u)&\leq& \sum_{(\vk{k},\vk{l})\in\mathcal{K}_3(u)}
\pk{\sup_{(\vk{x},\vk{y})\in D_{\vk{k}}(u)\times D_{\vk{l}}(u) }\overline{X}(\vk{x})+\overline{X}(\vk{x})>2(u-\mathbb{Q}_{11})} \nonumber\\
&\leq& \sum_{(\vk{k},\vk{l})\in\mathcal{K}_3(u)}
e^{-\frac{(2u-2\mathbb{Q}_{11}-\mathbb{Q}_{12})^2}{2(4-\delta)}}\nonumber\\
&\leq& \mathbb{Q}u^{2n} e^{-\frac{(2u-2\mathbb{Q}_{11}-\mathbb{Q}_{12})^2}{2(4-\delta)}}\nonumber\\
&=&o\LT(u^{2n-1}\Psi(u)\RT), \ u\rw\IF,
\EQN
where $\mathbb{Q}_{11}=\sup_{\vk{x}\in E}g(\vk{x})<\IF$ and $\mathbb{Q}_{12}=2\E{\sup_{\vk{x}\in E}\overline{X}(\vk{x})}<\IF$.\\
Inserting \eqref{upper}-\eqref{err3} into \eqref{up1} and \eqref{low1} yields that
\BQNY
\pk{\sup_{\vk{x}\in E(u)}(X(\vk{x})+g(\vk{x}))>u}\sim \LT(\Pi_{i=1}^{n-1}(T_i-S_i)\RT)\LT(\Pi_{i=1}^nc_i\RT) \int_{-\IF}^{\IF}e^{-bx^{2}+cx}dx u^{2n-1}\Psi(u), \ u\rw\IF,
\EQNY
which compared with \eqref{asym1} implies the final result.
\QED

\BL\label{lem0}
Let $X_{u,k}(\vk{x}), \ k\in K_u, u \in D(u)$  be a family of centered Gaussian fields \eE{with continuous sample paths} where  $D(u)=\Pi_{i=1}^n[0, \lambda_i u^{-2}]$ for some $\vk{\lambda}>\vk{0}$. Let further  $u_k,  k\in K_u$ be  given positive constants  satisfying
\BQN\label{uuk}
\lim_{u\rw\IF}\sup_{k\in K_u}\LT|\frac{u_k}{u}-1\RT|=0.
\EQN
\eE{If} $X_{u,k}$ has unit variance, and correlation function $r_k$ \eE{(not depending on u)} satisfying \eqref{r2} uniformly with respect to $k\in K_u$, then
\BQNY
\lim_{u\rw\IF}\sup_{k\in K_u}\LT|
\frac{\pk{\sup_{\vk{x}\in D(u)}X_{u,k}(\vk{x})>u_{k}}}{\Psi(u_{k})}
-\mathcal{H} _{\vk{c}}[\vk{0},\vk{\lambda}]
\RT|=0.
\EQNY
\EL
\prooflem{lem0} It follows along the same lines of \cite{Uniform2016}[Theorem 2.1].
\COM{
Conditioning on $\mathcal{A}_{u}(w,k):=\LT\{X_{u,k}(\vk{0})=
u_k-\frac{w}{u_k}\RT\}, w\in \R$, we have for all $u$ large enough,
\begin{align*}
&\frac{\pk{\sup_{\vk{x}\in D(u)}X_{u,k}(\vk{x})>u_{k}}}{\Psi(u_{k})}\\
&=\frac{1}{\sqrt{2\pi}u_k\Psi(u_{k})}\int_{\R}e^{-\frac{1}{2}\LT(u_k-\frac{w}{u_k}\RT)^2}
\pk{\sup_{\vk{x}\in [\vk{0},\vk{\lambda}]}X_{u,k}(u^{-2}\vk{x})>u_{k}\Bigl\lvert\mathcal{A}_{u}(w,k) }dw\\
&=\frac{e^{-\frac{u_k^2}{2}}}{\sqrt{2\pi}u_k\Psi(u_{k})}\int_{\R}
e^{w-\frac{w^2}{2u^2_k}}
\pk{\sup_{\vk{x}\in [\vk{0},\vk{\lambda}]}\mathcal{X}_u^w(\vk{x},k)>w}dw\\
&=:\frac{e^{-\frac{u_k^2}{2}}}{\sqrt{2\pi}u_k\Psi(u_{k})} I_{u,k},
\end{align*}
where $$\mathcal{X}_u^w(\vk{x},k)=u_k\LT(X_{u,k}(u^{-2}\vk{x})-u_k\RT)+w
\Bigl\lvert\mathcal{A}_{u}(w,k).$$
By \eqref{uuk}, we have
\BQNY
\lim_{u\rw\IF}\sup_{k\in K_u}\abs{\frac{e^{-\frac{u_k^2}{2}}}{\sqrt{2\pi}u_k\Psi(u_{k})}-1}=0.
\EQNY
Thus in order to establish the proof, it suffices to prove that
\BQN\label{IR}
\lim_{u\rw\IF}\sup_{k\in K_u}\LT|I_{u,k}- \Hn _{\vk{a}}[\vk{0},\vk{\lambda}]\RT|=0.
\EQN
For $W$ some positive constant, it follows that
\BQNY
&&\sup_{k\in K_u}\LT|I_{u,k}- \Hn _{\vk{a}}[\vk{0},\vk{\lambda}]\RT|\\
&&\leq \sup_{k\in K_u}\abs{\int_{-W}^{W}
\LT(e^{w-\frac{w^2}{2u^2_k}}
\pk{\sup_{\vk{x}\in [\vk{0},\vk{\lambda}]}\mathcal{X}_u^w(\vk{x},k)>w}-e^{w}
\pk{\sup_{\vk{x}\in [\vk{0},\vk{\lambda}]}\zeta(\vk{x})>w}\RT)dw}\\
&&\quad +\sup_{k\in K_u}\int_{(-\IF,-W]\cup[W,\IF)}e^{w-\frac{w^2}{2u^2_k}}
\pk{\sup_{\vk{x}\in [\vk{0},\vk{\lambda}]}\mathcal{X}_u^w(\vk{x},k)>w}dw\\
&&\quad +\int_{(-\IF,-W]\cup[W,\IF)}e^{w}
\pk{\sup_{\vk{x}\in [\vk{0},\vk{\lambda}]}\zeta(\vk{x})>w}dw\\
&&=:I_1+I_2+I_3,
\EQNY
where $\zeta(\vk{x})=Y(\vk{ax})-\Var(Y(\vk{ax}))$
with
\begin{align*}
Y(\vk{ax})=\sum_{i=1}^{n-1}B^i(a_ix_i)+B^{n}\LT(\sum_{i=2}^n a_ix_i\RT)
\end{align*}
where  $B^i(x),\ x\in\mathbb{R}, i=1,2\ldot   n$ are independent  standard Brownian motions.\\
Next, we give upper bounds for $I_i(u), i=1,2,3$.\\
\underline{Upper bound for $I_1(u)$}.
Direct calculations show that
\BQNY
\E{\mathcal{X}_u^w(\vk{x},k)}=-u_k^2(1-r_k(u^{-2}\vk{x}))+w(1-r_k(u^{-2}\vk{x}))
\EQNY
and
\BQNY
\Var\LT(\mathcal{X}_u^w(\vk{x},k)-\mathcal{X}_u^w(\vk{y},k)\RT)=
u_k^2\LT(\Var(X_{u,k}(\vk{x})-X_{u,k}(\vk{y}))
-\LT(r_k(u^{-2}\vk{x})-r_k(u^{-2}\vk{y})\RT)^2\RT).
\EQNY
By (\ref{r2}) and (\ref{uuk}), it follows  uniformly with respect to $\vk{x}\in[\vk{0},\vk{\lambda}], k\in K_u, w\in[-W,W]$ that
\BQN\label{XM12}
\E{\mathcal{X}_u^w(\vk{x},k)}\rw
-\sum_{i=1}^{n-1}\abs{a_i x_i}-\abs{\sum_{i=2}^{n}a_ix_i}
\EQN
 as $u\rw\IF$ and also for any $t,t'\in [S_1,S_2]$ uniformly with respect to $k\in K_u$, any $w_i\in \R$,
\BQN\label{XVar}
\Var\LT(\mathcal{X}_u^w(\vk{x},k)-\mathcal{X}_u^w(\vk{y},k)\RT)\rw 2\sum_{i=1}^{n-1}\abs{a_i (x_i-y_i)}+2\abs{\sum_{i=2}^{n}a_i(x_i-y_i)},
\EQN
as $u\rw\IF$.
Combination of (\ref{XM12}) and (\ref{XVar}) shows that the convergence of finite-dimensional distributions of
$$\mathcal{X}_u^w(\vk{x},k),\ \vk{x}\in[\vk{0},\vk{\lambda}]$$
is $\zeta(\vk{x}), \ \vk{x}\in[\vk{0},\vk{\lambda}]$.

Moreover, by (\ref{r2}) we have that there exists a constant $C>0$ such that for all $\vk{x}, \vk{y}\in[\vk{0},\vk{\lambda}]$ and  all large  $u$
\BQN\label{HC}
\sup_{k\in K_u}\Var\LT(\mathcal{X}_u^w(\vk{x},k)-\mathcal{X}_u^w(\vk{y},k)\RT)
\leq C\sum_{i=1}^{n}\abs{x_i-y_i},
\EQN
which combined with \eqref{XM12} implies that the family of distributions
$$\pk{\mathcal{X}_u^w(\vk{x},k)\in(\cdot)}$$
is uniformly tight with respect to $k\in K_u$ and $w$ in a compact set of $\R$. Consequently,
\BQNY
\{\mathcal{X}_u^w(\vk{x},k),\ \vk{x}\in[\vk{0},\vk{\lambda}]\} \quad \text{weakly converges to} \quad  \{\zeta(\vk{x}),\ \vk{x}\in[\vk{0},\vk{\lambda}]\}.
\EQNY
Let $$\mathbb{A}:=\left\{v: \pk{\sup_{\vk{x}\in[\vk{0},\vk{\lambda}]}(\zeta(\vk{x})-v)>0}  \text{is continuous at }  v\right\}.$$
Note that if $w\in \mathbb{A}$, then
 $$\pk{\sup_{\vk{x}\in[\vk{0},\vk{\lambda}]}(\zeta(\vk{x})-v)>x}$$
 is continuous with respect to $x$ at \eE{$x\not=0$}.
Hence by the continuity of supremum functional, we have that as
\BQNY
c_u(w):&=&\sup_{k\in K_u}\LT|\pk{\sup_{\vk{x}\in[\vk{0},\vk{\lambda}]}\mathcal{X}_u^w(\vk{x},k)>w}
-\pk{\sup_{\vk{x}\in[\vk{0},\vk{\lambda}]}\zeta(\vk{x})>w}\RT|\rw 0,\ \ u\rw\IF,
\EQNY
for $w\in \mathbb{A}.$
 Noting that $mes(\mathbb{A}^c)=0$ and by dominated convergence theorem, we have that
 \BQNY
 I_1(u)\leq e^{W}\int_{[-W,W] \cap \mathbb{A} }c_u(w)dw+2W e^{W}\sup_{\vk{w}\in [-W, W]}\left|1-e^{-\frac{w^2}{2u_k^2}}\right|\rw 0, \quad u\rw\IF.
 \EQNY

\underline{Upper bound for $I_2(u)$}.
Using (\ref{XM12}) and (\ref{XVar}), for some $\delta\in(0,1/2)$, $|w|>W$ with $W$ sufficiently large and all $u$ large we have
\BQNY
\sup_{k\in K_u,\vk{x}\in[0,\vk{\lambda}]}\E{\mathcal{X}_u^w(\vk{x},k)}
\leq\mathbb{C}_1+\delta|w|
\EQNY
and
\BQNY
\sup_{k\in K_u,\vk{x}\in[0,\vk{\lambda}]}\Var\LT(\mathcal{X}_u^w(\vk{x},k)\RT)
\leq \mathbb{C}_2.
\EQNY
by (\ref{HC}) and Theorem 8.1 of \cite{Pit96}
\BQNY
I_2(u)&=&\sup_{k\in K_u}\int_{(-\IF,-W]\cup[W,\IF)}e^{w-\frac{w^2}{2u^2_k}}
\pk{\sup_{\vk{x}\in [\vk{0},\vk{\lambda}]}\mathcal{X}_u^w(\vk{x},k)>w}dw\\
&\leq&e^{-W}+\sup_{k\in K_u}\int_{W}^{\IF)}e^{w}
\pk{\sup_{\vk{x}\in [\vk{0},\vk{\lambda}]}\LT(\mathcal{X}_u^w(\vk{x},k)
-\E{\mathcal{X}_u^w(\vk{x},k)}\RT)>(1-\delta)w-\mathbb{C}_1}dw\\
&\leq& e^{-W}+\int_W^{\IF}e^{w}\mathbb{C}_3{w}^2
\Psi\LT(\frac{(1-\delta)w-\mathbb{C}_1}{\mathbb{C}_2}\RT)d w\\
&=:&A_1(W)\rw 0,\ W\rw\IF,
\EQNY
\underline{Upper bound for $I_3(u)$}.
Borell-TIS inequality (see, e.g., \cite{AdlerTaylor}) implies that  $$I_3(u)\rw 0,\ u, W\rw\IF,$$
hence (\ref{IR}) follows.

\QED

\begin{remark}
Following \cite{EGRFRV2017} [Lemma 5.2], in \nelem{lem0} the results of the two dimensional cases holds.
\end{remark}
}

\BT \label{ThDK}
Assume that $f:[0,1]^n\to \mathbf{R}$ is a measurable bounded positive function. Let $$F(x_1\ldot x_n)=\int_0^{x_1}\dots \int_0^{x_n} f(t_1\ldot t_n)dt_1\dots dt_{n}$$ and assume that $F(1\ldot  1)=1$. Then the projection of the level set $\mathcal{D}=\{(x_1\ldot x_n): F(x_1\ldot x_n)=1/2\}$ on $\mathbb{R}^{n-1}$ is a Jordan measurable set of nonzero measure.

\ET
\prooftheo{ThDK}
The assumption on $f$ tells us that there is a positive constant $M$ so that
\BQN \label{bern}
0< f(x) \le  M
\EQN
 for  every point $x\in [0,1]^n$.
Since $f$ is defined on $[0,1]^n$ and positive, we can extend it on $\mathbb{R}^n$ subject to the condition \eqref{bern}. This can be done in an explicit way by using symmetries w.r.t. to all sides of the cube $[0,1]^n$. Namely, we say that two points $x=(x_1\ldot x_n),y=(y_1\ldot y_n)\in\mathbb{R}^n$ are equivalent if there are integers $i,k,m$ so that $x_i-y_i = 2k$ or $x_i+y_i=2m+1$. Now every $y\in\mathbb{R}^n$ has its unique representatives $x$ in $[0,1]^n$. Then we set $f(y)=f(x)$ for $y\not\in [0,1]^n$.

By the definition of $F$ we can find the unique global solution $$x_n=h(x_1\ldot x_{n-1})$$ of the equation $$F(x_1\ldot  x_{n-1},x_n)=\frac{1}{2}$$ for
$$x'=(x_1\ldot  x_{n-1})\in (0,\infty)^{n-1},$$ because $F$ is strictly increasing w.r.t all variables. Here $h$ is strictly decreasing w.r.t. all variables $x_1\ldot  x_{n-1}$. Moreover $h$ is continuous.

Further we are interested on the set $$L=\{x'=(x_1\ldot  x_{n-1}): x'\in [0,1]^{n-1}, F(x',h(x'))=1/2\}$$ and have to prove that it is Jordan measurable.

We have that $(x_1\ldot  x_{n-1})\in K$ if and only if $(x_1\ldot  x_{n-1})\in [0,1]^{n-1}$ and $0\le h(x_1\ldot  x_{n-1})\le 1$. The boundary of this set is contained in the union of the following two sets
$$A:=\partial [0,1]^n, \quad \text{ and  } B:=\{(x_1\ldot  x_{n-1})\in [0,1]^n: h(x_1\ldot  x_{n-1})=1\}.$$ By using again the previous procedure, we find that there is a  mapping $g:[0,1]^{n-2}\to [0,1]$ so that $$B=\{(x_1\ldot x_{n-2},g(x_1\ldot x_{n-2})): (x_1\ldot x_{n-2})\in[0,1]^{n-2}\}.$$ Let $T=\{(x_1\ldot  x_{n-2}, x_{n-1}): (x_1\ldot x_{n-2})\in[0,1]^{n-2}, 0\le x_{n-1}\le g(x_1\ldot x_{n-1})\}$. Since $g$ is continuous, it follows that $T$ is Jordan measurable and its boundary contains $B$. Namely $$\mu(T)=\int_{[0,1]^{n-2}} g(x_1\ldot  x_{n-2})dx_1\dots dx_{n-2}.$$ Thus $B$ has Jordan $n-1$ measure equal to zero. This implies that $\partial L$ has Lebesgue measure 0 (with respect to Lebesgues measure on $\R^{n-1}$ and in particular $L$ is measurable. Further we show that its measure is not zero. Since $F(1\ldot 1,1)=1$ and $F(0\ldot 0,0)=0$, it follows that there is $t\in(0,1)$ so that $F(t\ldot t,t)=1/2$. Then $h(t\ldot t)=t$, and thus by the continuity it follows that  $T=(t\ldot t)$ is an interior point of $L$, namely for $\epsilon>0$ small enough so that $t+\epsilon<1$,  there is a ball $B$ centered at $T$ and with a positive radius $\delta$ such that $h(B)\subset (t-\epsilon,t+\epsilon)$. Consequently,   $L$ has a positive Lebesgue measure establishing the proof. \QED

\section*{Acknowledgments} 
Thanks to Professor Enkelejd Hashorva who give many useful suggestions which greatly improve our manuscript. Thanks to Swiss National Science Foundation Grant no.  200021-175752.

	\bibliographystyle{ieeetr}
\bibliography{KSA}
\end{document}